\documentclass{amsart}

\usepackage[utf8]{inputenc}
\usepackage[T1]{fontenc}
\usepackage{lmodern}
\usepackage{amsmath}
\usepackage{amsthm}
\usepackage{amssymb}
\usepackage{units}
\usepackage{amscd}
\usepackage{geometry}
\usepackage{graphicx}
\usepackage{tikz}
\usepackage{tikz-cd}
\usepackage{pst-plot}
\usepackage{enumitem}
\usetikzlibrary{arrows,automata}
\usepackage{graphicx}
\usepackage{url}
\usepackage{pgfplots}
\usepackage{hyperref}

{

}
{\newtheorem{theo}{Theorem}[section]}
{\newtheorem{corollaire}[theo]{Corollary}}
{}
{}
{}
{}
{\theoremstyle{definition} \newtheorem{defin}[theo]{Definition}
							\newtheorem{construction}[theo]{Construction}
							\newtheorem{prop}[theo]{Proposition}
							\newtheorem{lemme}[theo]{Lemma}}
{\theoremstyle{remark} \newtheorem{remarque}[theo]{Remark}
						}

	\newcommand{\Id}{Id}	

\newcommand{\op}{op}

\newcommand{\sS}{\text{sSet}}
\newcommand{\Hom}{\text{Hom}}

\newcommand{\Top}{\text{Top}}
\newcommand{\Real}[1]{{}||#1||}
\newcommand{\RealP}[1]{{}\Real{#1}_P}
\newcommand{\Sing}{\text{Sing}}

\newcommand{\pr}{pr}
\newcommand{\colim}{\operatornamewithlimits{colim}}
\newcommand{\Ex}{\text{Ex}}
\newcommand{\sd}{\text{sd}}
\newcommand{\Exi}{\text{Ex}^{\infty}}
\newcommand{\lv}{\text{l.v}}
\newcommand{\Map}{\text{Map}}
\newcommand{\Fun}{\text{Fun}}

\newcommand{\nd}{\text{n.d.}}

\newcommand{\fil}[1]{(#1,\varphi_{#1})}

\newcommand{\C}{\textbf{C}}
\newcommand{\D}{\textbf{D}}

\newcommand{\N}{\mathbb{N}}
\newcommand{\Diag}{\text{Diag}}

\newcommand{\RealNP}[1]{\Real{#1}_{N(P)}}

\newcommand{\tr}{\text{tr}}

\newcommand{\Ho}{\text{Ho}}

\newcommand{\Sd}{\text{Sd}}
\newcommand{\ex}{\text{ex}}

\newcommand{\fib}{\text{fib}}
\newcommand{\Ver}{\text{Vert}}
\newcommand{\CW}{\text{CW}}
\newcommand{\filP}[1]{(#1,\varphi_P\circ\varphi_{#1})}
\newcommand{\lab}[1]{(#1,\lambda_{#1})}
\renewcommand{\Im}{\text{Im}}
\pgfplotsset{compat=1.16}

\title{A Stratified Kan-Quillen Equivalence}
\author{Sylvain Douteau}
\begin{document}
\begin{abstract}In this paper, we exhibit a Quillen equivalence between two model categories encoding the homotopy theory of stratified spaces : the model category of filtered simplicial sets, and that of filtered spaces. Additionally, we introduce a new class of filtered spaces, that of vertical filtered CW-complexes, providing a nice model for the homotopy category of stratified spaces
\end{abstract}
\maketitle

\tableofcontents
Let $\mathcal{C}$ be a category, in which we want to invert some class of morphisms $W\subset \mathcal{C}$. One can try to define directly the localized category $\mathcal{C}[W^{-1}]$, but even if it exists, this "homotopy category" is not the proper setting to do homotopical algebra. Indeed, many of the constructions that might be available in $\mathcal{C}$, such as colimits and limits, will not have well-behaved counterparts in $\mathcal{C}[W^{-1}]$. Instead, one would want some notion of derived (co)-limits, invariant with respect to maps in $W$. For this, one needs more than just the datum of the category $\mathcal{C}[W^{-1}]$. 

This is the problem that Quillen's model categories fix \cite{QuillenHomotopicalAlgebra}. In a model category, in addition to picking a class $W$ of weak-equivalences, one also picks two extra classes of maps, the cofibrations and fibrations. Those classes satisfy axioms that are inspired by the properties of cellular maps, Serre fibrations and weak homotopy equivalences. This leads to a natural construction of a homotopy relation between morphisms, and the realization of the localized category $\mathcal{C}[W^{-1}]$ as an explicit category, usually written $\Ho(\mathcal{C})$, whose objects are a subclass of the objects of $\mathcal{C}$ and whose maps are homotopy classes of maps in $\mathcal{C}$. In this setting, one can define derived functors, homotopy (co)-limits, or further localize with respect to additional morphisms. 

Moreover, Quillen defined equivalences between model categories, which are now known as Quillen equivalences. All previously mentioned constructions are invariant under Quillen equivalences. In particular, Quillen equivalent model categories have equivalent homotopy categories. For this reason, given a particular homotopy category $\Ho(\mathcal{C})$, it is very useful to find a model category $\mathcal{D}$ such that $\Ho(\mathcal{D})$ is equivalent to $\Ho(\mathcal{C})$, but it is even more important to understand the Quillen equivalence class of $\mathcal{D}$. Abstractly, a homotopy theory is a class of Quillen equivalent model categories.

Among model categories, simplicial model categories play an important role because they grant access to simplicial sets of morphisms. These enriched sets of morphisms are derived versions of the classical sets of morphisms. They are the main invariants of a homotopy theory, and provide a natural source for many homotopical invariants \cite{DwyerKanFunctionComplexes}.

%

The homotopy category corresponding to the classical homotopy theory of spaces can be described as the category of CW-complexes and homotopy classes of continuous maps. There are two important model categories corresponding to this homotopy category : the category of simplicial sets, $\sS$, together with the Kan model structure, and the category of topological spaces, $\Top$, together with the Quillen model structure. Those two model categories are related by a classical adjunction,
\begin{equation}\label{EquationKanQuillen}
\Real{-}\colon \sS\leftrightarrow \Top\colon \Sing
\end{equation}
and Quillen showed that it was a Quillen equivalence \cite{QuillenHomotopicalAlgebra}.
This means that both model categories encode the same homotopy theory, that of spaces.
In particular, any computation of homotopy limits and colimits, of homotopy localizations or of simplicial mapping spaces can be done in either model category.
This is especially important since the combinatorial nature of the category of simplicial sets grants access to a lot of natural constructions. For example, through the use of simplicial (abelian) groups, one can naturally define classifying spaces for groups, or loop spaces.

Nowadays, there is a growing interest in understanding the homotopy theory of stratified spaces, see \cite{AyalaFrancisTanaka,Haine,NandLal,Lukas}. This is in part motivated by the observation that classical topological objects can be better understood through the additional data of a stratification (see \cite{ArticleMoiWhitehead,Jansen} for some recent work on the subject).
Recall that stratified spaces are topological spaces, together with the data of some decomposition into strata, usually given as a continuous map toward a poset, $X\to P$. Historically, they were introduced to deal with manifolds with singularities, with the goal of extending invariants of manifolds to those objects. This is what is accomplished by Goresky and MacPherson's intersection homology \cite{IntersectionHomologyI}. Nowadays, finer invariants of those stratified spaces have been introduced, such as the $(\infty-)$category of exit-paths \cite{Treumann,WoolfFundamentalCategory,HigherAlgebra}. Those invariants are only well-behaved with respect to stratified maps and stratified homotopy equivalences, which leads to the study of the associated homotopy theory.

Following Quillen's strategy, the author of this paper defined two model categories for spaces stratified over some partially ordered set, $P$. One for the simplicial side, the model category of filtered simplicial sets $\sS_P$ \cite{ArticleMoi}, which is a simplicial model category, and one for the topological side, the model category of filtered spaces $\Top_P$ \cite{TopNP}. The former is similar to the Kan model structure on the category of simplicial sets : its cofibrations are the monomorphisms, its fibrations are characterized by lifting conditions against horns, and its weak-equivalences are detected by suitably defined homotopy groups. The latter is similar to the Quillen model structure on topological spaces. All objects are fibrant, weak-equivalences are defined as maps inducing isomorphisms on all (filtered) homotopy groups and fibrations are characterized by lifting conditions similar to those defining Serre fibrations. Furthermore, the adjunction \eqref{EquationKanQuillen} has an equivalent in the stratified setting :
\begin{equation}\label{EquationStratifiedKanQuillen}
\RealP{-}\colon \sS_P\leftrightarrow \Top_P\colon \Sing_P
\end{equation}
Given the similarity between the classical and the stratified context, one might expect this adjunction to be a Quillen equivalence, but it is in fact not even a Quillen adjunction since $\Sing_P$ does not preserve fibrant objects (see \cite[Example 4.14]{ArticleMoi}). Nevertheless, it is known that using a suitably defined filtered subdivision, one recovers a Quillen adjunction \cite[Corollaire 8.2.4]{TheseMoi},
\begin{equation}\label{EquationQEAMontrer}
\RealP{\sd_P(-)}\colon \sS_P\leftrightarrow\Top_P\colon\Ex_P\Sing_P
\end{equation}
In this paper, we prove that \textbf{the Quillen adjunction \eqref{EquationQEAMontrer} is a Quillen equivalence.}
This result fits into a larger corpus of work, investigating the homotopy theory of stratified spaces and a stratified version of the homotopy hypothesis (see \cite{AyalaFrancisTanaka,Haine,NandLal,Lukas}). This theorem unites two different points of view on the homotopy theory of stratified spaces. A topological point of view, carried by the category of filtered spaces, $\Top_P$, in which pseudo-manifolds and maps between them live. And a combinatorial point of view, given by the simplicial model category $\sS_P$. The latter is the natural context in which to interpret the higher invariants of stratified spaces, such as the $\infty$-category of exit-paths. It should also be noted that this results does not follow from a straight-forward adaptation of Quillen's original proof of the equivalence between $\Top$ and $\sS$. Instead, it relies on a comparison with a third model category of diagrams of simplicial sets.

The proof of this result is broken down into two parts. One on the simplicial side, the other on the topological side. Along the way, we prove that several useful adjunctions are also Quillen equivalences. Notably, we fill in a gap in the proof of \cite[Theorem 2.15]{TopNP}.
In addition, we also introduce vertical filtered CW-complexes, a nice class of stratified spaces, which provides a model for the homotopy category of stratified spaces (see Corollary \ref{CorollaireCategorieHomotopieCW}).

\section{Preliminaries and sketch of proof}
In this section, we recall the needed definitions and theorems from \cite{ArticleMoi} and \cite{TopNP}, and give a sketch of the proof of the main theorem.

\subsection{Filtered simplicial sets and diagrams}
\begin{defin}
Let $P$ be a partially ordered set. A filtered simplicial set over $P$ is the data of a simplicial set $X$, together with  a simplicial map toward the nerve of $P$, $\varphi_X\colon X\to N(P)$.
A filtered map between filtered simplicial sets $f\colon \fil{X}\to\fil{Y}$ is a simplicial map $f\colon X\to Y$ such that the following triangle commutes
\begin{equation*}
\begin{tikzcd}
X\arrow{rr}{f}
\arrow[swap]{dr}{\varphi_X}
&& Y
\arrow{dl}{\varphi_Y}
\\
&N(P)
\end{tikzcd}
\end{equation*}
The category of filtered simplicial sets over $P$, $\sS_P$, is the category of filtered simplicial sets over $P$ and filtered maps.
\end{defin}

\begin{defin}
A filtered simplex is the data of a simplicial map $\varphi\colon \Delta^n\to N(P)$. Such a filtered simplex will usually be denoted $\Delta^{\varphi}$ instead of $(\Delta^n,\varphi)$. Alternatively, by identifying $\Delta^{\varphi}$ with the image of $\varphi$ in $N(P)$, we will write $\Delta^{\varphi}=[p_0,\dots,p_n]$. A non-degenerate filtered simplex is the data of an injective simplicial map $\varphi\colon \Delta^n\to N(P)$. Alternatively, filtered simplices can be seen as simplices of the simplicial set $N(P)$. From this point of view, one can associate to each filtered simplex $\varphi\colon \Delta^n\to N(P)$, the unique non-degenerate simplex of $N(P)$, $\Delta^{\bar{\varphi}}$ such that $\Delta^{\varphi}$ is a degeneracy of $\Delta^{\bar{\varphi}}$. We write $\Delta(P)$ and $R(P)$ for the full sub-categories of $\sS_P$, generated by the filtered simplices and the non-degenerate filtered simplices, respectively.
\end{defin}

\begin{remarque}
There is an equivalence of categories $\sS_P\simeq \text{Fun}(\Delta(P)^{\op},\text{Set})$.
\end{remarque}

\begin{defin}
Let $\Diag_P$ be the category of functors :
\begin{equation*}
\Diag_P=\Fun(R(P)^{\op},\sS)
\end{equation*}
Furthermore, we consider $\Diag_P$ with its projective model structure, in which weak-equivalences and fibrations are determined level-wise. It is a cofibrantly generated model category (see \cite[proposition 2.2]{TopNP} for the generating sets of (trivial) cofibrations).
\end{defin}

\begin{defin}
Define the bi-functor $-\otimes-\colon \sS\times\sS_P\to\sS_P$ as
\begin{equation*}
K\otimes\fil{X}\mapsto (K\times X,\varphi_X\circ\pr_X)
\end{equation*}
If $\fil{X}\in \sS_P$ is fixed, the functor $-\otimes\fil{X}\colon \sS\to \sS_P$ admits a right adjoint $\Map(\fil{X},-)\colon \sS_P\to \sS$ defined on objects as
\begin{equation*}
\Map(\fil{X},\fil{Y})_n=\Hom_{\sS_P}(\Delta^n\otimes\fil{X},\fil{Y})
\end{equation*}
\end{defin}

Recall that $\sS_P$ equipped with the functors : $-\otimes -\colon \sS\times\sS_P\to \sS_P$, and $\Map(-,-)\colon \sS_P^{\op}\times\sS_P\to \sS$, is a simplicial category (see \cite[section 3]{ArticleMoi}). 

\begin{defin}
Define $D$ as the functor
\begin{align*}
D\colon \sS_P&\to\Diag_P\\
\fil{X}&\mapsto \left\{\begin{array}{ccl}
R(P)^{\op}&\to &\sS\\
\Delta^{\varphi}&\mapsto &\Map(\Delta^{\varphi},\fil{X})
\end{array}\right.
\end{align*}
where $\Map$ is the simplicial $\Hom$ in $\sS_P$.
The functor $D$ admits a left adjoint, $C\colon\Diag_P\to\sS_P$, defined as follows. Let $\mathcal{C}$ be the full subcategory of $R(P)^{\op}\times R(P)$ whose objects are the pairs $(\Delta^{\varphi},\Delta^{\psi})$ such that $\Delta^{\psi}\subset\Delta^{\varphi}$. To any functor in $\Diag_P$, $F\colon R(P)^{\op}\to \sS$, associate the following functor
\begin{align*}
F\otimes R(P)\colon \mathcal{C}&\to \sS_P\\
(\Delta^{\varphi},\Delta^{\psi})&\mapsto F(\Delta^{\varphi})\otimes\Delta^{\psi}
\end{align*}
Then, define $C(F)$ as the colimit :
\begin{equation*}
C(F)=\colim_{\mathcal{C}}F\otimes R(P)
\end{equation*}
\end{defin}


One can then define a model structure on $\sS_P$ (\cite[Theorem 3.19]{ArticleMoi})
\begin{theo}
The category $\sS_P$ admits a simplicial combinatorial model structure where :
\begin{itemize}
\item the cofibrations are the monomorphisms,
\item the fibrations are defined by lifting conditions against "admissible horn"
\item the weak-equivalence are the maps $f\colon X\to Y$ such that $D(f^{\fib})\colon D(X^{\fib})\to D(Y^{\fib})$ is a weak-equivalence in $\Diag_P$, where $(-)^{\fib}$ is any fibrant replacement functor.
\end{itemize}
\end{theo}

We will also make heavy use of two endofunctors of $\sS_P$, $\sd_P$ and $\Ex_P$, for which we quickly recall the definitions. See \cite[section 2.1]{ArticleMoi} for complete definitions.
\begin{defin}
Let $P$ be a partially ordered set. Let $\sd_P(N(P))\subset \sd(N(P))\times N(P)$ be the sub-simplicial set containing the following simplices :
\begin{equation*}
\sd_P(N(P))=\{((\Delta^{\varphi_0},\dots,\Delta^{\varphi_k}),\Delta^{\psi})\ |\ \Delta^{\varphi_0}\subset\dots\subset\Delta^{\varphi_k}\subset N(P),\ \psi\colon \Delta^k\to N(P),\  \psi(\Delta^k)\subset \Delta^{\varphi_0}\}
\end{equation*}
Alternatively, under the identification $\Delta^{\psi}=[p_0,\dots,p_k]$ we can describe $\sd_P(N(P))$ as 
\begin{equation*}
\{((\Delta^{\varphi_0},p_0),\dots,(\Delta^{\varphi_k},p_k))\ |\ \Delta^{\varphi_0}\subset\dots\subset\Delta^{\varphi_k}\subset N(P),\ p_0\leq\dots\leq p_k, p_i\in\Delta^{\varphi_0}, \forall\ 0\leq i\leq k\}
\end{equation*}
This simplicial set comes equiped with the two projections : $\pr_1\colon \sd_P(N(P))\to \sd(N(P))$ and $\pr_2\colon \sd_P(N(P))\to N(P)$.

Let $\fil{X}$ be a filtered simplicial set. Its filtered subdivision $\sd_P\fil{X}$ has for underlying simplicial set the following pull-back :
\begin{equation*}
\begin{tikzcd}
sd_P\fil{X}
\arrow{r}
\arrow[swap]{d}{\sd_P(\varphi_X)}
&\sd(X)
\arrow{d}{\sd(\varphi_X)}
\\
\sd_P(N(P))
\arrow{r}{\pr_1}
&\sd(N(P))
\end{tikzcd}
\end{equation*}
The filtration on $\sd_P\fil{X}$ is then given by the composition 
\begin{equation*}
\begin{tikzcd}
\sd_P\fil{X}\arrow{r}{\sd_P(\varphi_X)} &\sd_P(N(P))\arrow{r}{\pr_2} &N(P)
\end{tikzcd}
\end{equation*}
This definition extends to a well-defined functor
\begin{equation*}
\sd_P\colon \sS_P\to \sS_P
\end{equation*}
Furthermore, there is a natural transformation $\lv_P\colon \sd_P\to \Id$. The functor $\sd_P$ admits a right adjoint, $\Ex_P$, which comes with a natural transformation $\beta\colon\Id\to \Ex_P$.
\end{defin}

\subsection{Filtered and strongly filtered spaces}

Throughout this paper, $\Top$ denotes the category of $\Delta$-generated spaces (see \cite{Dugger}, and \cite{ConvenientCategory}). All topological spaces are assumed to be $\Delta$-generated.

\begin{defin}
Let $P$ be a partially ordered set. A filtered space over $P$ is the data of 
\begin{itemize}
\item a space, $X$,
\item a continuous map $\varphi_X\colon X\to P$, where $P$ is given the Alexandrov topology.
\end{itemize}
A filtered map $f\colon\fil{X}\to\fil{Y}$ is a continuous map $f\colon X\to Y$, such that the following triangle commutes
\begin{equation*}
\begin{tikzcd}
X
\arrow{rr}{f}
\arrow[swap]{dr}{\varphi_X}
&&Y
\arrow{dl}{\varphi_Y}
\\
&P
\end{tikzcd}
\end{equation*}
The category of filtered spaces and filtered maps is denoted $\Top_P$.
\end{defin}

\begin{defin}
A strongly filtered space over $P$ is the data of 
\begin{itemize}
\item A topological space $X$,
\item A continuous map $\varphi_X\colon X\to \Real{N(P)}$, where $\Real{N(P)}$ is the realization of the simplicial set $N(P)$.
\end{itemize}
A strongly filtered map $f\colon \fil{X}\to\fil{Y}$ is a continuous map $f\colon X\to Y$ such that the following commutative triangle commutes :
\begin{equation*}
\begin{tikzcd}
X
\arrow[swap]{dr}{\varphi_X}
\arrow{rr}{f}
&&Y
\arrow{dl}{\varphi_{Y}}
\\
&\Real{N(P)}
\end{tikzcd}
\end{equation*}
The category of strongly filtered spaces and strongly filtered maps is denoted $\Top_{N(P)}$.
\end{defin}

\begin{defin}
Let $\varphi_P\colon \Real{N(P)}\to P$ be the continuous map, defined by, if $t$ is a point in the interior of $\Delta^{\psi}=[p_0,\dots,p_k]\subset N(P)$, $\varphi_P(t)= p_k$. Composition of the filtration with $\varphi_P\colon \Real{N(P)}\to P$ induces a functor
\begin{align*}
\varphi_P\circ-\colon \Top_{N(P)}&\to \Top_P\\
(X,\varphi_X\colon X\to \Real{N(P)})&\mapsto (X,\varphi_P\circ\varphi_X\colon X\to P)
\end{align*}
A right adjoint to this functor, $-\times_P\Real{N(P)}\colon \Top_P\to \Top_{N(P)}$, is defined on objects as the following pull-back :
\begin{equation*}
\begin{tikzcd}
\fil{Y}\times_P\Real{N(P)}
\arrow{r}
\arrow{d}
&\Real{N(P)}
\arrow{d}
\\
\fil{Y}
\arrow{r}
&P
\end{tikzcd}
\end{equation*}
The strong filtration is then given by projection on the second factor.
\end{defin}

We then have the expected adjunctions relating (strongly) filtered spaces and filtered simplicial sets.

\begin{defin}
Let $\fil{X}$ be a filtered simplicial set. Define $\RealNP{\fil{X}}\in \Top_{N(P)}$ as 
\begin{equation*}
\RealNP{\fil{X}}=(\Real{X},\Real{\varphi_X}\colon \Real{X}\to\Real{N(P)})
\end{equation*}
By defining $\RealNP{-}$ on maps as the usual realisation, we get a functor 
\begin{equation*}
\RealNP{-}\colon \sS_P\to \Top_{N(P)}
\end{equation*}
Define $\RealP{-}\colon \sS_P\to \Top_P$ as the composition $(\varphi_P\circ -) \RealNP{-}$
\end{defin}

\begin{prop}
The functor $\RealNP{-}$ admits a right adjoint, $\Sing_{N(P)}\colon \Top_{N(P)}\to \sS_P$, such that, for any strongly filtered space, $\fil{X}$, the following diagram is a pullback square
\begin{equation*}
\begin{tikzcd}
\Sing_{N(P)}\fil{X}
\arrow[hookrightarrow]{r}
\arrow{d}
&\Sing(X)
\arrow{d}{\Sing(\varphi_X)}
\\
N(P)
\arrow[hookrightarrow]{r}
&\Sing(\Real{N(P)})
\end{tikzcd}
\end{equation*}
The composition $\Sing_P= \Sing_{N(P)} (-\times_P\Real{N(P)})$ is then a right adjoint to $\RealP{-}$.
\end{prop}

\begin{defin}
Define the functors $\C$ as the composition $\C = \RealNP{C(-)} \colon \Diag_P\to \Top_{N(P)}$, and $\D$ as the composition $\D = D \ \Sing_{N(P)}\colon \Top_{N(P)}\to \Diag_P$.
\end{defin}

\subsection{Proof of the main theorem and outline of the paper}

In this paper, we prove the following theorem :

\begin{theo}
The adjunction
\begin{equation*}
\RealP{\sd_P(-)}\colon \sS_P\leftrightarrow \Top_P\colon \Ex_P\Sing_P
\end{equation*}
is a Quillen-equivalence.
\end{theo}

\begin{proof}
The proof of the main theorem can be summarized by the following diagram.
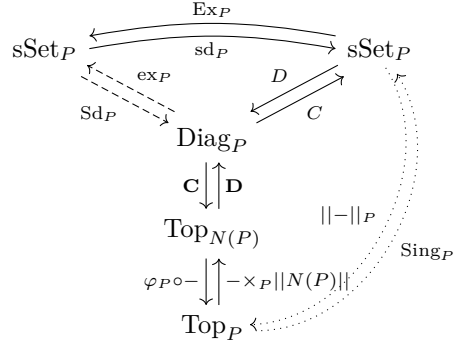
\begin{equation*}
\begin{tikzcd}
\sS_P
\arrow[shift right,swap = 1, bend left = 10]{rr}{\sd_P}
\arrow[shift right,swap = 1,dashed]{dr}{\Sd_P}
&
&\sS_P
\arrow[shift right,swap = 1, bend right = 10]{ll}{\Ex_P}
\arrow[shift right,swap = 1]{dl}{D}
\arrow[shift right, swap, bend left = 60, dotted]{dddl}{\RealP{-}}
\\
&\Diag_P
\arrow[shift right,swap = 1, dashed]{ul}{\ex_P}
\arrow[shift right,swap = 1]{ur}{C}
\arrow[shift right,swap = 1]{d}{\textbf{C}}
\\
&\Top_{N(P)}
\arrow[shift right,swap = 1]{u}{\textbf{D}}
\arrow[shift right,swap = 1]{d}{\varphi_P\circ-}
\\
&\Top_P
\arrow[shift right,swap = 1]{u}{-\times_P \Real{N(P)}}
\arrow[shift right, swap, dotted, bend right = 60]{uuur}{\Sing_P}
\end{tikzcd}
\end{equation*}

In section \ref{SubsectionFactoringSd}, we will construct a pair of adjoint functor $\Sd_P\colon \sS_P\leftrightarrow \Diag_P\colon \ex_P$ such that $C\circ \Sd_P\simeq \sd_P$, and $\ex_P\circ D\simeq \Ex_P$ (Proposition \ref{PropSdPFactorisesdP}). This implies that the functor $\Ex_P\Sing_P$ is isomorphic to $\ex_PD\Sing_P$. By \cite[Remark 4.16]{ArticleMoi}, the latter is isomorphic to $\ex_P(\D(-\times_P\Real{N(P)}))$. In particular, it suffices to show that all three of the adjunctions : $(\Sd_P,\ex_P)$, $(\C,\D)$, and $(\varphi_P\circ-,-\times_P\Real{N(P)})$ are Quillen equivalences.

\begin{itemize}
\item The adjunction $(\Sd_P,\ex_P)$ is defined and studied in section \ref{SectionSd}, and by Theorem \ref{TheoSdPQuillenAdjunction}, it is a Quillen equivalence. We also show in this section that the adjoint pairs $(C,D)$ and $(\sd_P,\Ex_P)$ are Quillen equivalences (Theorem \ref{TheoCDQuillenEquivalence} and Proposition \ref{QEsdEx})
\item The adjunction $(\C,\D)$ is a Quillen equivalence by \cite[Theorem 2.12]{TopNP}.
\item The adjunction $(\varphi_P\circ-,-\times_P\Real{N(P)})$ is studied in section \ref{SectionTopNP}. The proof of Theorem 2.15 in \cite{TopNP}, stating that it is a Quillen equivalence was missing arguments. To fill the gap, we introduce the notion of vertical filtered CW-complexes in subsection \ref{SubsectionCW}, and use them to prove that it is indeed a Quillen equivalence (Theorem \ref{TheoremeAMontrer}).
\end{itemize}
\end{proof}
%
%

\section{Filtered simplicial sets and subdivisions}
\label{SectionSd}
In this section, we investigate the adjunctions between the category of filtered simplicial sets $\sS_P$ and the category of diagrams, $\Diag_P$. In subsection \ref{subsectionCD}, we first prove that the adjoint pair $(C,D)$ is a Quillen equivalence (Theorem \ref{TheoCDQuillenEquivalence}). Then, in subsection \ref{SubsectionFactoringSd}, we describe a pair of adjoint functors $(\Sd_P,\ex_P)$, factoring the adjoint pair $(\sd_P,\Ex_P)$ through $(C,D)$ (Propositions \ref{PropSdPFunctor} and \ref{PropSdPFactorisesdP}), and we show that it is a Quillen equivalence (Theorem \ref{TheoSdPQuillenAdjunction}).

\subsection{Equivalence between diagrams and filtered simplicial sets}
\label{subsectionCD}
In this subsection, we prove the following theorem.

\begin{theo}\label{TheoCDQuillenEquivalence}
The adjoint pair $(C,D)$ is a Quillen equivalence.
\end{theo}

We will make use the following lemma, from \cite{ArticleMoi}.

\begin{lemme}[{\cite[Proposition 3.6, Proposition 3.8]{ArticleMoi}}]
The pair $(C,D)$ is a Quillen adjunction. Furthermore, $D$ reflects weak-equivalences between fibrant objects.
\end{lemme}

\begin{proof}[Proof of Theorem \ref{TheoCDQuillenEquivalence}]
By \cite[Corollary 1.3.16(c)]{Hovey}, since $D$ reflects weak-equivalences between fibrant objects, it is enough to show that for any cofibrant diagram $F$, the map
\begin{equation*}
F\to D((C(F))^{\fib})
\end{equation*}
is a weak-equivalence, where $(-)^{\fib}$ is any fibrant replacement functor. Since weak-equivalences in $\Diag_P$ are defined levelwise, one needs to check for all $\Delta^{\varphi}\in R(P)$ that the map
\begin{equation*}
F(\Delta^{\varphi})\to \Map(\Delta^{\varphi},(C(F))^{\fib})
\end{equation*}
is a weak-equivalence of simplicial sets.
Since $F$ is cofibrant, by lemma \ref{LemmeExPVarphiPart}, it is enough to check that the map 
\begin{equation*}
F(\Delta^{\varphi})\to \Map(\Delta^{\varphi},(F(\Delta^{\varphi})\otimes\Delta^{\varphi})^{\fib})
\end{equation*}
Is a weak-equivalence. Now, let $F(\Delta^{\varphi})\to K$ be some fibrant replacement of $F(\Delta^{\varphi})$ (say, $K=\Exi(F(\Delta^{\varphi}))$), and consider the following commutative diagram
\begin{equation*}
\begin{tikzcd}
F(\Delta^{\varphi})
\arrow{r}
\arrow{d}
&
\Map(\Delta^{\varphi},(F(\Delta^{\varphi})\otimes\Delta^{\varphi})^{\fib})
\arrow{d}
\\
K
\arrow{r}
& \Map(\Delta^{\varphi},(K\otimes\Delta^{\varphi})^{\fib})
\end{tikzcd}
\end{equation*}
By construction, the left arrow is a trivial cofibration.
 Furthermore, since $-\otimes\Delta^{\varphi}$ preserves trivial cofibrations,
  so is $F(\Delta^{\varphi})\otimes\Delta^{\varphi}\to K\otimes\Delta^{\varphi}$.
   In addition, $(-)^{\fib}$ also preserves weak-equivalences, and takes value in fibrant objects.
    Since $D$ preserves weak-equivalences between fibrant objects, the right arrow is also a weak-equivalence.
     By two out of three, it is then sufficient to show that the bottom arrow is a weak-equivalence.
      First, notice that $K\otimes\Delta^{\varphi}\to (K\otimes\Delta^{\varphi})^{\fib}$ is a weak-equivalence between
       fibrant objects. Which means that the map 
\begin{equation*}
\Map(\Delta^{\varphi},K\otimes\Delta^{\varphi})\to \Map(\Delta^{\varphi},(K\otimes\Delta^{\varphi})^{\fib})
\end{equation*} 
       is a weak-equivalence. Since it also factors the bottom arrow, it is sufficient to show that the following map is a weak-equivalence.
\begin{equation*}
K\to \Map(\Delta^{\varphi},K\otimes\Delta^{\varphi})
\end{equation*}
But the codomain of this map is isomorphic to $\Map(\Delta^{|\varphi|},K)$ which is homotopy equivalent to $K$.
\end{proof}

\begin{lemme}\label{LemmeExPVarphiPart}
Let $F$ be a cofibrant object in $\Diag_P$, and $\Delta^{\varphi}$ be a non-degenerated simplex in $R(P)$. The map $F(\Delta^{\varphi})\otimes\Delta^{\varphi}\to C(F)$ induces a weak-equivalence
\begin{equation*}
\Map(\Delta^{\varphi},(F(\Delta^{\varphi})\otimes\Delta^{\varphi})^{\fib})\to \Map(\Delta^{\varphi},(C(F))^{\fib})
\end{equation*}
for any fibrant replacement functor $(-)^{\fib}$.
\end{lemme}

The proof of Lemma \ref{LemmeExPVarphiPart} will rely on the following definition and lemmas. 

\begin{defin}
Let $\fil{X}$ be a filtered simplicial set, and $\Delta^{\varphi}\in R(P)$ a non-degenerate filtered simplex. We say that a filtered simplex $\sigma\colon \Delta^{\psi}\to\fil{X}$ is of $\varphi$-type, if $\Delta^{\psi}$ is a degeneracy of $\Delta^{\varphi}$. We define the $\varphi$-part of $X$, $X^{\varphi}\subset X$, as the subsimplicial set of $X$ generated by the simplices of $\varphi$-type. $(X^{\varphi},(\varphi_X)_{|X^{\varphi}})$ is a filtered sub-simplicial set of $\fil{X}$, that we will denote $\fil{X}^{\varphi}$. This defines a functor
\begin{equation*}
(-)^{\varphi}\colon \sS_P\to \sS_P
\end{equation*}
\end{defin}

\begin{lemme}\label{LemmeMapVarphiPart}
Let $\fil{X}$ be a filtered simplicial set and $\Delta^{\varphi}$ a non-degenerate filtered simplex. Then, the inclusion $\fil{X}^{\varphi}\to\fil{X}$ induces an isomorphism
\begin{equation*}
\Map(\Delta^{\varphi},\fil{X}^{\varphi})\to \Map(\Delta^{\varphi},\fil{X})
\end{equation*}
\end{lemme}

\begin{proof}
Let $\fil{X}$ be a filtered simplicial set and $\Delta^{\varphi}\in R(P)$ a non-degenerate filtered simplex. A simplex in $\Map(\Delta^{\varphi},\fil{X})$ is a filtered simplicial map $\sigma\colon \Delta^n\otimes\Delta^{\varphi}\to\fil{X}$.
Passing to the $\varphi$-part, we get
\begin{equation*}
(\sigma)^{\varphi}\colon (\Delta^n\otimes\Delta^{\varphi})^{\varphi}\to\fil{X}^{\varphi}
\end{equation*}
But any simplex in $\Delta^n\otimes\Delta^{\varphi}$ is a face of a simplex of $\varphi$-type, which means that $(\Delta^n\otimes\Delta^{\varphi})^{\varphi}=\Delta^n\otimes\Delta^{\varphi}$. In particular, all simplices of $\Map(\Delta^{\varphi},\fil{X})$ factor through $\fil{X}^{\varphi}$.
\end{proof}

\begin{lemme}\label{LemmeVarphiPartLeftQuillen}
The functor $(-)^{\varphi}\colon \sS_P\to\sS_P$ preserves weak-equivalences and fibrant objects.
\end{lemme}

\begin{proof}
We will prove both assertion separately. First, note that $(-)^{\varphi}$ preserves all colimits, which means it is a left adjoint. We will show that it is a left Quillen functor. This will prove the first claim, since a left Quillen functor preserves weak-equivalences between cofibrant objects, and all objects are cofibrant in $\sS_P$. Since $\sS_P$ is cofibrantly generated, it is enough to show that the generating (trivial) cofibrations are sent to (trivial) cofibrations. Since, $(-)^{\varphi}$ clearly preserves monomorphism, it is true for cofibration. Let $\Lambda^{\psi}_k\to \Delta^{\psi}$ be an admissible horn inclusion. $\Delta^{\psi}$ must be a degenerate simplex of $N(P)$, let $\Delta^{\bar{\psi}}\in R(P)$ be the non-degenerate simplex such that $\Delta^{\psi}$ is a degeneracy of $\Delta^{\bar{\psi}}$. Consider the following cases :
\begin{itemize}
\item If $\Delta^{\varphi}\not\subset\Delta^{\bar{\psi}}$ then no simplex in $\Delta^{\psi}$ is of $\varphi$-type. In particular, we have
\begin{equation*}
\left(\Lambda^{\psi}_k\to\Delta^{\psi}\right)^{\varphi}=\emptyset\to \emptyset
\end{equation*}
which is a weak-equivalence.
\item If $\Delta^{\varphi}\subsetneq\Delta^{\bar{\psi}}$, then, the top dimensional simplex of $\Delta^{\psi}$ is not of $\varphi$-type. Furthermore, since $\Lambda^{\psi}_k\to\Delta^{\psi}$ is admissible, $d_k(\Delta^{\psi})$ is of $\bar{\psi}$-type, which means it is also not of $\varphi$-type. In particular, $(\Lambda^{\psi}_k)^{\varphi}=(\Delta^{\psi})^{\varphi}$, and the generating trivial cofibration is sent to an isomorphism.
\item If $\Delta^{\varphi}=\Delta^{\bar{\psi}}$, then $(\Delta^{\psi})^{\varphi}=\Delta^{\psi}$. Furthermore, one has 
\begin{equation*}
(\Lambda^{\psi}_k)^{\varphi}=\coprod\limits_{\substack{j\not = k\\ d_j\Delta^{\psi} \text{ of $\varphi$-type}}}d_j\Delta^{\psi}
\end{equation*}
This is a generalized horn in the sense of definition \ref{DefinGeneralizedHorn}, and since $\Lambda^{\psi}_k$ is an admissible horn, there exists $\epsilon=\pm 1$ such that $\psi(k)=\psi(k+\epsilon)$. We conclude by lemma \ref{LemmeGeneralizedHorn}.
\end{itemize}

For the second part of Lemma \ref{LemmeVarphiPartLeftQuillen}, let $\fil{X}$ be a fibrant object in $\sS_P$. We need to show that $\fil{X}^{\varphi}$ is fibrant. By \cite[Lemma 1.17]{ArticleMoi}, it is enough to show that the map $\fil{X}^{\varphi}\to N(P)$ admits the right lifting property against all maps of the form 
\begin{equation*}
\Delta^1\otimes\partial(\Delta^{\psi})\cup \{\epsilon\}\otimes\Delta^{\psi}\to\Delta^1\otimes\Delta^{\psi},
\end{equation*}
Where $\Delta^{\psi}$ is any filtered simplex, and $\{\epsilon\}$ is one of the two faces of $\Delta^1$. Consider such a lifting problem :
\begin{equation}\label{EquationLiftingAlternateHorn}
\begin{tikzcd}
\Delta^1\otimes\partial(\Delta^{\psi})\cup \{\epsilon\}\otimes\Delta^{\psi}
\arrow{d}
\arrow{r}{f}
&\fil{X}^{\varphi}
\arrow{d}
\\
\Delta^1\otimes\Delta^{\psi}
\arrow{r}
\arrow[dotted]{ur}{g}
&N(P)
\end{tikzcd}
\end{equation}
Once again, we distinguish three cases, depending on which type of simplex $\Delta^{\psi}$ is. We will write $\Delta^{\bar{\psi}}$ for the non-degenerate simplex corresponding to $\Delta^{\psi}$.
\begin{itemize}
\item If $\Delta^{\bar{\psi}}\not\subset\Delta^{\varphi}$, then there exists no lifting problem of the form \eqref{EquationLiftingAlternateHorn}. Indeed, by definition of $\fil{X}^{\varphi}$, $\Hom(\Delta^{\psi},\fil{X}^{\varphi})=\emptyset$.
\item If $\Delta^{\bar{\psi}}=\Delta^{\varphi}$, consider the composition of the top map in \eqref{EquationLiftingAlternateHorn} with the inclusion $\fil{X}^{\varphi}\to \fil{X}$. Since, $\fil{X}$ is fibrant, there exists a lift $\Delta^1\otimes\Delta^{\psi}\to \fil{X}$. But $\Delta^{\psi}$ is of $\varphi$-type, and so are the top dimensional simplices of $\Delta^1\otimes\Delta^{\psi}$, so this lift lands in $\fil{X}^{\varphi}$, providing a lift in \eqref{EquationLiftingAlternateHorn}.
\item If $\Delta^{\bar{\psi}}\subsetneq\Delta^{\varphi}$, consider the restriction of $f$ to $\{\epsilon\}\otimes\Delta^{\psi}$, $f_{|\{\epsilon\}\otimes\Delta^{\psi}}\colon \Delta^{\psi}\to \fil{X}^{\varphi}$. By definition of $\fil{X}^{\varphi}$, there must exist some simplex $\sigma\colon \Delta^{\mu}\to \fil{X}$ such that $\Delta^{\mu}$ is of $\varphi$-type, and the restriction of $f$ is a face of $\sigma$. Furthermore, $\sigma$ lands in $\fil{X}^{\varphi}$. Consider the modified lifting problem :
\begin{equation*}
\begin{tikzcd}[column sep = huge]
\Delta^1\otimes\partial(\Delta^{\psi})\cup \{\epsilon\}\otimes\Delta^{\mu}
\arrow{d}
\arrow{r}{f\cup \left(\{\epsilon\}\otimes\sigma\right)}
&\fil{X}^{\varphi}
\arrow{r}
&\fil{X}
\arrow{dl}
\\
\Delta^1\otimes\Delta^{\mu}
\arrow{r}
\arrow[dashed,near start]{urr}{\widetilde{g}}
&N(P)
\arrow[crossing over, from=u]
\end{tikzcd}
\end{equation*}
The leftmost arrow is a cofibration, and since its domain and codomain retract to $\{\epsilon\}\otimes\Delta^{\mu}$, it is also a stratified homotopy equivalence, which means it is a trivial cofibration. Since $\fil{X}$ is fibrant, there must exist a lift, $\widetilde{g}\colon \Delta^1\otimes\Delta^{\mu}\to\fil{X}$. Since $\Delta^{\mu}$ is of $\varphi$-type, so are the top dimensional simplices of $\Delta^1\otimes\Delta^{\mu}$, which means $\widetilde{g}$ lands in $\fil{X}^{\varphi}$. Restricting $\widetilde{g}$ to $\Delta^1\otimes\Delta^{\psi}\subset \Delta^1\otimes \Delta^{\mu}$ gives a lift of \eqref{EquationLiftingAlternateHorn}, $g\colon \Delta^1\otimes\Delta^{\psi}\to \fil{X}^{\varphi}$.
\end{itemize}
\end{proof}

\begin{defin}\label{DefinGeneralizedHorn}
Let $\psi\colon \Delta^n\to N(P)$ be a filtered simplex, and $S$ a non-empty, proper subset of $\{0,\dots,n\}$.
Define the generalized horn $\Lambda^\psi_S$ as the sub-simplicial set of $\Delta^{\psi}$ generated by the following faces
\begin{equation*}
\Lambda^\psi_S=\coprod_{j\not\in S}d_j(\Delta^{\psi})
\end{equation*}
\end{defin}
\begin{lemme}\label{LemmeGeneralizedHorn}
Let $\Lambda^{\psi}_S$ be a generalized horn. If there exists $j$ and $\epsilon =\pm 1$ such that $j\in S$, $j+\epsilon \not \in S$ and $\psi(j)=\psi(j+\epsilon)$, then the inclusion $\Lambda^{\psi}_S\to \Delta^{\psi}$ is a trivial cofibration.
\end{lemme}

\begin{proof}
We will prove this statement by induction on the cardinality of $S$.
Let $\psi\colon \Delta^n\to N(P)$ be a filtered simplex and $S\subset \{0,\dots,n\}$ be a non-empty proper subset. Since $S$ is non-empty, we have $|S|\geq 1$. If $|S|=1$, $S=\{s\}$ and the generalized horn, $\Lambda^{\psi}_S$, is just the horn $\Lambda^{\psi}_S$. By hypothesis, there exists $\epsilon =\pm 1$ such that $\psi(s)=\psi(s+\epsilon)$, which means the horn $\Lambda^{\psi}_s$ is an admissible horn, and the maps $\Lambda^{\psi}_s\to \Delta^{\psi}$ is a trivial cofibration.

Now, assume $|S|>1$. By hypothesis, there exists $j\in S$ satisfying the hypothesis. Pick $a\not= j\in S$, and let $A=S\setminus \{a\}$. The map  $\Lambda^{\psi}_S\to \Delta^{\psi}$ factors through the inclusion $\Lambda^{\psi}_S\to \Lambda^{\psi}_A$ which fits into the following pushout square
\begin{equation*}
\begin{tikzcd}
d_a(\Delta^{\psi})
\cap \Lambda^{\psi}_S
\arrow{r}
\arrow{d}
&\Lambda^{\psi}_S
\arrow{d}
\\
d_a(\Delta^{\psi})
\arrow{r}
&\Lambda^{\psi}_A
\end{tikzcd}
\end{equation*}
Since $|A|=|S|-1$, it is enough to show that $\Lambda^{\psi}_S\to\Lambda^{\psi}_A$ is a trivial cofibration. Since trivial cofibrations are stable under pushout, we need to show that the leftmost map is a trivial cofibration. Notice that $d_a(\Delta^{\psi})$ is isomorphic to the following filtered simplex 
\begin{align*}
\psi'\colon \Delta^{n-1}&\to N(P)\\
i&\mapsto \left\{\begin{array}{cl}
\psi(i) &\text{ if $i<a$}\\
\psi(i+1) &\text{ if $i\geq a$}
\end{array}\right.
\end{align*}
Furthermore, under this isomorphism, $d_a(\Delta^{\psi})\cap \Lambda^{\psi}_S$ is isomorphic to $\Lambda^{\psi'}_B$, where $B=(D^a)^{-1}(A)$ (with $D^a\colon \{0,\dots,n-1\}\to \{0,\dots,n\}$ the strictly increasing map omitting $a$). By construction $|B|=|A|<|S|$. Furthermore, by setting $j'=(D^a)^{-1}(j)\in B$, one has $j'+\epsilon=(D^a)^{-1}(j+\epsilon)\in B$ and $\psi'(j')=\psi(j)=\psi(j+\epsilon)=\psi'(j'+\epsilon)$. 
In particular, $B$ satisfy the hypothesis of Lemma \ref{LemmeGeneralizedHorn}, so by the induction hypothesis $\Lambda^{\psi'}_B\to \Delta^{\psi'}$ is a trivial cofibration.
\end{proof}

\begin{lemme}\label{LemmeVarphiPartColim}
Let $F$ be a cofibrant object in $\Diag_P$, and $\Delta^{\varphi}\in R(P)$ a non-degenerate filtered simplex. The map, 
\begin{equation*}
F(\Delta^{\varphi})\otimes\Delta^{\varphi}\to C(F)
\end{equation*}
is injective, and its image is $(C(F))^{\varphi}$.
\end{lemme}
\begin{proof}
Since $F$ is a cofibrant object in $\Diag_P$, by \cite[Lemma 2.13]{TopNP}, for any non-degenerate filtered simplex $\Delta^{\varphi}\in R(P)$, the map $i_{\varphi}\colon F(\Delta^{\varphi}\otimes\Delta^{\varphi})\to C(F)$ is a monomorphism. Furthermore, by construction, 
\begin{equation*}
C(F)=\bigcup\limits_{\Delta^{\varphi}\in R(P)} i_{\varphi}(F(\Delta^{\varphi})\otimes\Delta^{\varphi})
\end{equation*}
Now, given a simplex of $\varphi$-type, 
$\sigma\colon \Delta^{\psi}\to C(F)$, this simplex must land in some 
$i_{\mu}(F(\Delta^{\mu})\otimes\Delta^{\mu})$. 
Furthermore, for such a map to exist, one must have $\Delta^{\varphi}\subset\Delta^{\mu}$. 
But then, if $\Delta^{\mu}\not=\Delta^{\varphi}$, 
$\sigma$ must land in 
$i_{\mu}(F(\Delta^{\mu})\otimes\Delta^{\varphi})$, which is identified in 
$C(F)$ with a subsimplicial set of 
$i_{\varphi}(F(\Delta^{\varphi})\otimes\Delta^{\varphi})$. In particular, all simplices of 
$\varphi$-type in $C(F)$ are in 
$i_{\varphi}(F(\Delta^{\varphi})\otimes\Delta^{\varphi})$. 
Since, conversely, $F(\Delta^{\varphi})\otimes\Delta^{\varphi}$ is generated by simplices of 
$\varphi$-type, we have 
$(C(F))^{\varphi}=i_{\varphi}(F(\Delta^{\varphi})\otimes\Delta^{\varphi})$, which concludes the proof.
\end{proof}
\begin{proof}[Proof of Lemma \ref{LemmeExPVarphiPart}]
Let $F$ be a cofibrant object in $\Diag_P$, and $\Delta^{\varphi}\in R(P)$ a non-degenerate filtered simplex.
Consider the following commutative diagram
\begin{equation*}
\begin{tikzcd}
(C(F))^{\varphi}
\arrow{r}
\arrow{d}
&C(F)
\arrow{d}
\\
((C(F))^{\varphi})^{\fib}
\arrow{r}
&(C(F))^{\fib}
\end{tikzcd}
\end{equation*}
The vertical map are weak-equivalences. Since $(-)^{\varphi}$ preserve weak-equivalences by lemma \ref{LemmeVarphiPartLeftQuillen}, this is also true in the following diagram
\begin{equation*}
\begin{tikzcd}
(C(F))^{\varphi}
\arrow{r}
\arrow{d}
&(C(F))^{\varphi}
\arrow{d}
\\
(((C(F))^{\varphi})^{\fib})^{\varphi}
\arrow{r}
&((C(F))^{\fib})^{\varphi}
\end{tikzcd}
\end{equation*}
Since the top map is the identity, by two out of three, the bottom map must also be a weak-equivalence. Furthermore, by lemma \ref{LemmeVarphiPartLeftQuillen}, its domain and codomain are fibrant object. Since $D$ is right Quillen, the map
\begin{equation*}
\Map(\Delta^{\varphi},(((C(F))^{\varphi})^{\fib})^{\varphi})
\to
\Map(\Delta^{\varphi},((C(F))^{\fib})^{\varphi})
\end{equation*}
must be a weak-equivalence. But then, by lemma \ref{LemmeMapVarphiPart}, one can rewrite the domain and codomain as follows :
\begin{equation*}
\Map(\Delta^{\varphi},((C(F))^{\varphi})^{\fib})
\to
\Map(\Delta^{\varphi},(C(F))^{\fib})
\end{equation*}
Finally, by lemma \ref{LemmeVarphiPartColim}, we have a weak-equivalence
\begin{equation*}
\Map(\Delta^{\varphi},(F(\Delta^{\varphi}\otimes\Delta^{\varphi}))^{\fib})
\to
\Map(\Delta^{\varphi},(C(F))^{\fib})
\end{equation*}
\end{proof}

\subsection{factoring the subdivision}
\label{SubsectionFactoringSd}

\begin{defin}
Let $\Delta^{\varphi}\in R(P)$ be a non-degenerate filtered simplex. Define the simplicial set $\Sd_P(N(P))(\Delta^{\varphi})$ as the image of the map
\begin{equation*}
\left(\sd_P(N(P))\right)^{\varphi}\to \sd(N(P))
\end{equation*}
Let $\fil{X}$ be a filtered simplicial set. Define the simplicial set $\Sd_P\fil{X}(\Delta^{\varphi})$ as the following pullback
\begin{equation*}
\begin{tikzcd}
\Sd_P\fil{X}(\Delta^{\varphi})
\arrow{r}
\arrow{d}
&\sd(X)
\arrow{d}{\sd(\varphi_X)}
\\
\Sd_P(N(P))(\Delta^{\varphi})
\arrow{r}
&\sd(N(P))
\end{tikzcd}
\end{equation*}
\end{defin}

\begin{prop}\label{PropSdPFunctor}
Let $\Delta^{\psi}\subset\Delta^{\varphi}$ be two non-degenerate filtered simplices and $\fil{X}$ a filtered simplicial set. There is an inclusion $\Sd_P\fil{X}(\Delta^{\varphi})\subset\Sd_P\fil{X}(\Delta^{\psi})$. In particular, $\Sd_P$ is a functor
\begin{equation*}
\Sd_P\colon \sS_P\to\Diag_P
\end{equation*}
furthermore, $\Sd_P$ has a right adjoint $\ex_P\colon \Diag_P\to \sS_P$.
\end{prop}

\begin{proof}
First, note that $\Sd_P(N(P))(\Delta^{\varphi})=\{(\sigma_0,\dots,\sigma_n)\ |\ \Delta^{\varphi}\subset \sigma_0\}\subset \sd(N(P))$. Indeed, given such a simplex of $\sd(N(P))$, if $\Delta^{\varphi}=[p_0,\dots,p_k]$, then $((\sigma_0,p_0),\dots,(\sigma_0,p_k),(\sigma_1,p_k),\dots,(\sigma_n,p_k))$ is a simplex in $\left(\sd_P(N(P))\right)^{\varphi}$  whose image in $\sd(N(P))$ is $(\sigma_0,\dots,\sigma_n)$.
In particular, if $\Delta^{\psi}\subset\Delta^{\varphi}$ are two non-degenerate filtered simplices, there is an inclusion of subsimplicial sets of $\sd(N(P))$ :
\begin{equation*}
\Sd_P(N(P))(\Delta^{\varphi})\subset \Sd_P(N(P))(\Delta^{\psi})
\end{equation*}
Next, let $\fil{X}$ be a filtered simplicial set, and consider the following commutative diagram :
\begin{equation*}
\begin{tikzcd}
&\Sd_P\fil{X}(\Delta^{\psi})
\arrow{ddd}
\arrow{dr}
\\
&&\sd(X)
\arrow{ddd}
\\
\Sd_P\fil{X}(\Delta^{\varphi})
\arrow{ddd}
\arrow[hookrightarrow,dashed]{uur}
\arrow[crossing over]{urr}
\\
&\Sd_P(N(P))(\Delta^{\psi})
\arrow{dr}
\\
&&\sd(N(P))
\\
\Sd_P(N(P))(\Delta^{\varphi})
\arrow{urr}
\arrow[hookrightarrow]{uur}
\end{tikzcd}
\end{equation*}
By construction, the squares on the front and on the back are pullback squares. This means that there is an unique map $\Sd_P\fil{X}(\Delta^{\varphi})\to \Sd_P\fil{X}(\Delta^{\psi})$ making the diagram commute. Furthermore, this map make the leftmost square into a pullback square. In particular, it must be a monomorphism. We have shown that $\Sd_P\fil{X}\colon R(P)^{\op}\to\sS$ is a functor. Furthermore, the functoriality of $\sd$, and of pullbacks, implies that $\Sd_P\colon \sS_P\to \Diag_P$ is a functor. Since pullback commute with colimits, $\Sd_P$ is a colimit-preserving functor between presheaf category, and it admits a right adjoint, that we will write $\ex_P\colon\Diag_P\to \sS_P$.
\end{proof}

\begin{prop}\label{PropSdPFactorisesdP}
There is an isomorphism of functor $C\circ\Sd_P\simeq \sd_P$
\end{prop}

\begin{proof}
Let $\fil{X}$ be a filtered simplicial set. By lemma \ref{LemmeDecompositionVarphiPart}, $\sd_P\fil{X}$ is naturally isomorphic to the colimit of the functor
\begin{align*}
G\colon \mathcal{C}&\to \sS_P\\
(\Delta^{\varphi},\Delta^{\psi})&\mapsto \left((\sd_P\fil{X})^{\varphi}\right)^{\psi}
\end{align*}
By lemma \ref{LemmeSdProduit}, the functor $G$ is isomorphic to $\Sd_P\fil{X}\otimes R(P)$. Since by definition, $C(\Sd_P\fil{X})= \colim\Sd_P\fil{X}\otimes R(P)$, we have a natural isomorphism
\begin{equation*}
C(\Sd_P\fil{X})\simeq \sd_P\fil{X}
\end{equation*}
\end{proof}

\begin{lemme}\label{LemmeDecompositionVarphiPart}
Let $\fil{X}$ be a filtered simplicial set. Define the functor
\begin{align*}
G\colon \mathcal{C}&\to \sS_P\\
(\Delta^{\varphi},\Delta^{\psi})&\mapsto (\fil{X}^{\varphi})^{\psi}
\end{align*}
Then, $\fil{X}\simeq \colim_{\mathcal{C}} G$.
\end{lemme}

\begin{proof}
Let us first prove that $G\colon \mathcal{C}\to \sS_P$ is well-defined. If $(\Delta^{\varphi_1},\Delta^{\psi_1}), (\Delta^{\varphi_2},\Delta^{\psi_2})\in \mathcal{C}$ are such that $\Delta^{\varphi_2}\subset\Delta^{\varphi_1}$ and $\Delta^{\psi_1}\subset\Delta^{\psi_2}$ (that is, if there is a map in $\mathcal{C}$ $(\Delta^{\varphi_1},\Delta^{\psi_1})\to (\Delta^{\varphi_2},\Delta^{\psi_2})$), we have $(\fil{X}^{\varphi_1})^{\psi_1}\subset (\fil{X}^{\varphi_2})^{\psi_2}$. Indeed, we have the following sequence of inclusion (using lemma \ref{LemmeCalculsVarphiPart})
\begin{align*}
\fil{X}^{\varphi_1}&\subset \fil{X}\\
(\fil{X}^{\varphi_1})^{\varphi_2}&\subset \fil{X}^{\varphi_2}\\
(\fil{X}^{\varphi_1})^{\psi_2}&\subset(\fil{X}^{\varphi_2})^{\psi_2}\\
\fil{X}^{\varphi_1})^{\psi_1}&\subset (\fil{X}^{\varphi_2})^{\psi_2}
\end{align*}
So the functor $G$ is well-defined. Then, by lemma \ref{LemmeCalculsVarphiPart}, if $\Delta^{\psi}\subset\Delta^{\varphi}$ are non-degenerate simplices, $(\fil{X}^{\varphi})^{\psi}= \fil{X}^{\varphi}\cap \fil{X}^{\psi}$. In particular, $\colim G$ corresponds to the gluing of the $\fil{X}^{\varphi}$ along their intersection. In particular, it is a subspace of $\fil{X}$, and it is enough to show that any simplex in $X$ is in some $\fil{X}^{\varphi}$. But if $\sigma\colon \Delta^{\psi}\to \fil{X}$ is a filtered simplex, then $\sigma$ lands in the $\bar{\psi}$-part of $\fil{X}$, where $\Delta^{\psi}$ is a degeneracy of the non-degenerate simplex $\Delta^{\bar{\psi}}$, which concludes the proof.
\end{proof}

\begin{lemme}\label{LemmeCalculsVarphiPart}
Let $\fil{X}$ be a filtered simplicial set, and $\Delta^{\varphi_1},\Delta^{\varphi_2}$ and $\Delta^{\varphi_3}$ non-degenerate filtered simplices. We have the following equalities and inclusion.
\begin{align*}
(\fil{X}^{\varphi_1})^{\varphi_1}&=\fil{X}^{\varphi_1}&\\
(\fil{X}^{\varphi_1})^{\varphi_2}&=\emptyset &\text{ if $\Delta^{\varphi_2}\not\subset \Delta^{\varphi_1}$}\\
(\fil{X}^{\varphi_1})^{\varphi_2}&=\fil{X}^{\varphi_1}\cap\fil{X}^{\varphi_2} &\text{ if $\Delta^{\varphi_2}\subset \Delta^{\varphi_1}$}\\
(\fil{X}^{\varphi_1})^{\varphi_3}&\subset (\fil{X}^{\varphi_1})^{\varphi_2} &\text{ if $\Delta^{\varphi_3}\subset \Delta^{\varphi_2}$}\\
((\fil{X}^{\varphi_1})^{\varphi_2})^{\varphi_3}&=(\fil{X}^{\varphi_1})^{\varphi_3} &\text{if $\Delta^{\varphi_3}\subset\Delta^{\varphi_2}\subset\Delta^{\varphi_1}$}
\end{align*}
\end{lemme}

\begin{proof}
The first equality follows from the definition of $\fil{X}^{\varphi_1}$. 
Then, if $\Delta^{\varphi_2}\not\subset\Delta^{\varphi_1}$, 
$\fil{X}^{\varphi_1}$ 
contains no filtered simplices of type $\varphi_2$, which implies the second equality. Now, let $\Delta^{\varphi_2}\subset\Delta^{\varphi_1}$ be two non-degenerate filtered simplices. There are two obvious inclusions $(\fil{X}^{\varphi_1})^{\varphi_2}\subset \fil{X}^{\varphi_1}$ and $\fil{X}^{\varphi_2}$, which gives half of the third equality. Conversely, let $\sigma$ be a simplex in $\fil{X}^{\varphi_1}\cap \fil{X}^{\varphi_2}$. Since $\sigma\in \fil{X}^{\varphi_1}$, there must be a simplex of type $\varphi_1$, $\tau\colon \Delta^{\psi}\to \fil{X}$, such that $\sigma$ is some face of $\tau$. But since $\sigma$ is also in $\fil{X}^{\varphi_2}$, $\sigma$ must be equal to some face of $\tau$ of type $\varphi_3$, with $\Delta^{\varphi_3}\subset\Delta^{\varphi_2}\subset\Delta^{\varphi_1}$. In particular, $\sigma$ must be a subface of a face of $\tau$ of type $\varphi_2$, and so $\sigma\in (\fil{X}^{\varphi_1})^{\varphi_2}$. This same reasoning implies the fourth assertion. The fifth follows from applying the third equality to $\fil{X}^{\varphi_1}$, then using the fourth assertion. 
\begin{align*}
((\fil{X}^{\varphi_1})^{\varphi_2})^{\varphi_3}&=(\fil{X}^{\varphi_1})^{\varphi_2}\cap(\fil{X}^{\varphi_1})^{\varphi_3}\\
&=(\fil{X}^{\varphi_1})^{\varphi_3}
\end{align*}
\end{proof}

\begin{lemme}\label{LemmeSdProduit}
Let $\fil{X}$ be a filtered simplicial set and $\Delta^{\varphi}$ a non-degenerate simplex. There is a natural isomorphism, 
\begin{equation*}
\Sd_P\fil{X}(\Delta^{\varphi})\otimes \Delta^{\varphi}\simeq(\sd_P\fil{X})^{\varphi}
\end{equation*}
\end{lemme}

\begin{proof}
Let us first prove the claim for $\fil{X}=N(P)$. The filtered simplicial set $(\sd_P(N(P)))^{\varphi}$ contains exactly the simplices of the form
\begin{equation}\label{EquationSimplicesVarphiPartsdP}
((\sigma_0,q_0),\dots,(\sigma_k,q_k)), \ \Delta^{\varphi}\subset\sigma_0, \  \{q_0,\dots,q_k\}\subset \{p_0,\dots,p_n\}
\end{equation}
with $\Delta^{\varphi}=[p_0,\dots,p_n]$. Indeed, let $((\sigma_0,q_0),\dots,(\sigma_k,q_k))$ be such a simplex. It is enough to construct a simplex of $\sd_P(N(P))$, $((\tau_0,r_0),\dots,(\tau_l,r_l))$ such that for all $0\leq i\leq k$ there exists some $0\leq j\leq l$ such that $(\sigma_i,q_i)=(\tau_j,r_j)$, and such that $\{r_0,\dots,r_l\}=\{p_0,\dots,p_n\}$. But one can obtain such a simplex by adding vertices of the form $(\sigma_k, p_i)$ for every $k$ such that $q_k<p_i$ and $q_{k+1}>p_i$, so all simplices of the form \eqref{EquationSimplicesVarphiPartsdP} are faces of simplices of $\varphi$-type. Conversely, a simplex of $\varphi$-type of $\sd_P(N(P))$  must be of the form $((\sigma_0,q_0),\dots,(\sigma_k,q_k))$, with $\{q_0,\dots,q_k\}=\{p_0,\dots,p_n\}$. Since it is a simplex of $\sd_P(N(P))$ one must also have $q_i\in \sigma_0$ for all $0\leq i\leq k$, and so $\Delta^{\varphi}\subset \sigma_0$. In particular, all faces of a simplex of $\varphi$-type of $\sd_P(N(P))$ satisfy the conditions of \eqref{EquationSimplicesVarphiPartsdP}. But \eqref{EquationSimplicesVarphiPartsdP} also describes the simplices of the product $\Sd_P(N(P))\otimes\Delta^{\varphi}$ (see the description of $\Sd_P(N(P))$ in the proof of Proposition \ref{PropSdPFunctor}), which concludes the proof for the case $\fil{X}=N(P)$.

Let $\fil{X}$ be a filtered simplicial set, consider the following commutative diagram
\begin{equation*}
\begin{tikzcd}
(\sd_P\fil{X})^{\varphi}
\arrow{r}
\arrow{d}
&\Sd_P\fil{X}(\Delta^{\varphi})
\arrow{r}
\arrow{d}
&\sd(X)
\arrow{d}
\\
(\sd_P(N(P)))^{\varphi}
\arrow{r}
&\Sd_P(N(P))(\Delta^{\varphi})
\arrow{r}
&\sd(N(P))
\end{tikzcd}
\end{equation*}
By definition of $\Sd_P$, the rightmost square is a pullback square. Furthermore, it is enough to show that the leftmost square is a pullback square since the isomorphism $(\sd_P(N(P)))^{\varphi}\simeq \Sd_P(N(P))(\Delta^{\varphi})\otimes\Delta^{\varphi}$ will then induce the desired isomorphism. We will show that the outer square is a pullback square. This is true by definition of $\sd_P$ if one omits the exponent $(-)^{\varphi}$, so it is enough to show that $\sd_P(\varphi_X)^{-1}((\sd_P(N(P)))^{\varphi})=(\sd_P\fil{X})^{\varphi}$. The inclusion $(\sd_P\fil{X})^{\varphi}\subset \sd_P(\varphi_X)^{-1}((\sd_P(N(P)))^{\varphi})$ follows from the fonctoriality of $(-)^{\varphi}$. On the other hand, fix some simplex in $\sd_P(\varphi_X)^{-1}((\sd_P(N(P)))^{\varphi})$. It must be of the form 
\begin{equation*}
(\sigma,((\sigma_0,q_0),\dots,(\sigma_k,q_k))), \ \sigma\colon \Delta^{\psi}\to \fil{X}, \sigma_0\subset\dots\subset\sigma_k\subset \Delta^{\psi},\ q_i\in \sigma_0, \ 0\leq i\leq k
\end{equation*}
with $\sigma_0$ of type $\mu$, with $\Delta^{\varphi}\subset\Delta^{\mu}$, and $\{q_0,\dots,q_k\}\subset\{p_0,\dots,p_n\}$. Proceeding as before, we can obtain this simplex as the face of some simplex of $\varphi$-type by adding vertices of the form $(\sigma_k, p_i)$ for the $p_i$ not equal to any of the $q_j$. In particular, any simplex in $\sd_P(\varphi_X)^{-1}((\sd_P(N(P)))^{\varphi})$ is the face of a simplex of $\varphi$-type of $\sd_P\fil{X}$. This prove that all square in the previous commutative diagram are pullback squares, which completes the proof.
\end{proof}

\begin{theo}\label{TheoSdPQuillenAdjunction}
The adjoint pair $(\Sd_P,\ex_P)$ is a Quillen equivalence.
\end{theo}

\begin{proof}
Let us first prove that it is a Quillen adjunction.
It is enough to show that $\Sd_P$ sends the generating (trivial) cofibrations of $\sS_P$ to (trivial) cofibrations in $\Diag_P$. We will start with the cofibrations. Let $\partial(\Delta^{\psi})\to\Delta^{\psi}$ be the inclusion of the boundary of some filtered simplex $\Delta^{\psi}$. For all $\Delta^{\varphi}\in R(P)$, define $F^{\varphi}$ to be the diagram defined by
\begin{align*}
R(P)^{\op}&\to \sS\\
\Delta^{\mu}&\mapsto \left\{\begin{array}{cl}
\Sd_P(\Delta^{\psi})(\Delta^{\varphi})&\text{ if $\Delta^{\mu}\subset\Delta^{\varphi}$}\\
\emptyset &\text{ if $\Delta^{\mu}\not\subset\Delta^{\varphi}$}
\end{array}\right.
\end{align*}
 We will decompose the map $\Sd_P(\partial(\Delta^{\psi}))\to \Sd_P(\Delta^{\psi})$ as a sequence of cofibrations as follows. First, define $F^{n+1}=\Sd_P(\partial(\Delta^{\psi}))$, where $n$ is the dimension of the image of $\Delta^{\psi}$ in $N(P)$. Then, for any $\Delta^{\varphi}\in R(P)$, write $|\varphi|$ for its dimension. For $0\leq k\leq n$, define $F^k$ as the following pushout 
\begin{equation*}
\begin{tikzcd}
\bigcup\limits_{|\varphi|=k}F^{\varphi}\cap F^{k+1}
\arrow{r}
\arrow{d}
&F^{k+1}
\arrow{d}
\\
\bigcup\limits_{|\varphi|=k}F^{\varphi}
\arrow{r}
&F^k
\end{tikzcd}
\end{equation*}
We will show that $F^0\simeq \Sd_P(\Delta^{\psi})$. Since by construction, the latter is covered by the $F^{\varphi}$, it is enough to show that for any $\Delta^{\varphi}\not=\Delta^{\mu}$ such that $|\varphi|=|\mu|=k$, $F^{\varphi}\cap F^{\mu}\subset F^{k+1}$. First, if $F^{\varphi}\cap F^{\mu}\not= \emptyset$, one must have $\Delta^{\varphi}\cap\Delta^{\mu}\not= \emptyset$, write this simplex $\Delta^{\nu}$. Furthermore, the intersection must be generated by $F^{\varphi}(\Delta^{\nu})\cap F^{\mu}(\Delta^{\nu})$. By the proof of proposition \ref{PropSdPFunctor}, one sees that all simplices of this intersection must be of the form $(\sigma_0,\dots,\sigma_l)$, with $\Delta^{\varphi}\subset\psi(\sigma_0)$ and $\Delta^{\mu}\subset\psi(\sigma_0)$ Since $\Delta^{\varphi}\not=\Delta^{\mu}$, there must exist some non-degenerate filtered simplex, $\Delta^{\epsilon}$, containing both $\Delta^{\varphi}$ and $\Delta^{\mu}$, and such that $\Delta^{\epsilon}\subset\sigma_0$. But then, $|\epsilon|>k$, which means the simplex is in $F^{k+1}$. In particular, $F^0\simeq \Sd_P(\Delta^{\psi})$, and the map $\Sd_P(\partial(\Delta^{\psi}))\to\Sd_P(\Delta^{\psi})$ can be obtained by successive pushout along maps of the form $F^{\varphi}\cap F^{k+1}\to F^{\varphi}$. In particular, it is enough to show that those are cofibrations. But the former diagram is of the form
\begin{align*}
R(P)^{\op}&\to\sS\\
\Delta^{\mu}&\mapsto \left\{\begin{array}{cl}
F^{\varphi}\cap F^{k+1}(\Delta^{\varphi}) &\text{ if $\Delta^{\mu}\subset\Delta^{\varphi}$}\\
\emptyset \text{ if $\Delta^{\mu}\not\subset\Delta^{\varphi}$}
\end{array}\right.
\end{align*}
where $F^{\varphi}\cap F^{k+1}(\Delta^{\varphi})$ is the subset of $F^{\varphi}(\Delta^{\varphi})=\Sd_P(\Delta^{\psi})(\Delta^{\varphi})$, defined by 
\begin{equation*}
F^{\varphi}\cap F^{k+1}(\Delta^{\varphi})=\{(\sigma_0,\dots,\sigma_l)\ |\ \text{$(\Delta^{\varphi'}\subset\sigma_0$, and $\Delta^{\varphi}\subsetneq\Delta^{\varphi'}$) or ($\Delta^{\varphi}\subset\sigma_0$ and $\sigma_l\not=\Delta^{\psi}$)}\}.
\end{equation*}
In particular, the maps $F^{\varphi}\cap F^{k+1}\to F^{\varphi}$ are cofibrations (see \cite[Proposition 2.2]{TopNP}), which concludes the proof in the case of cofibrations.

Now, consider a generating trivial cofibration in $\sS_P$, an admissible horn inclusion $\Lambda^{\psi}_k\to\Delta^{\psi}$. 
We already know that its image by $\Sd_P$ is a cofibration, (between cofibrant objects) so it is enough to show that it is also a weak equivalence. 
From \cite[Proposition 2.9]{ArticleMoi}, we know that $\sd_P(\Lambda^{\psi}_k)\to \sd_P(\Delta^{\psi})$ is a weak-equivalence in $\sS_P$. Furthermore, by Proposition \ref{PropSdPFactorisesdP} it is the image of the map $\Sd_P(\Lambda^{\psi}_k)\to\Sd_P(\Delta^{\psi})$ by the functor $C$. By Theorem \ref{TheoCDQuillenEquivalence}, $C$ is a left adjoint in a Quillen equivalence, so it reflects weak-equivalences between cofibrant objects, which implies that $\Sd_P(\Lambda^{\psi}_k)\to\Sd_P(\Delta^{\psi})$ is a weak-equivalence.

It then follows from Theorem \ref{TheoCDQuillenEquivalence} and Proposition \ref{QEsdEx}, and the two out of three properties of Quillen equivalences among Quillen adjunction that $(\Sd_P,\ex_P)$ is a Quillen equivalence.

\end{proof}

\begin{prop}\label{QEsdEx}
The adjoint pair $(\sd_P, \Ex_P)$ is a Quillen equivalence.
\end{prop}

\begin{proof}
By definition, $\sd_P$ preserve cofibrations. Furthermore, by \cite[Proposition 2.9]{ArticleMoi}, it also preserves trivial cofibrations, so it is a left Quillen functor. Now, let $f\colon \fil{X}\to\fil{Y}$ be some map in $\sS_P$. Consider the commutative diagram
\begin{equation*}
\begin{tikzcd}
\fil{X}
\arrow{r}{f}
\arrow[swap]{d}{\lv_P}
&\fil{Y}
\arrow{d}{\lv_P}
\\
\sd_P\fil{X}
\arrow{r}{\sd_P(f)}
&\sd_P\fil{Y}
\end{tikzcd}
\end{equation*}

By \cite[Lemma A.3]{ArticleMoi}, the vertical maps are weak-equivalences, which means by the two out of three property that $f$ is a weak-equivalence if and only if $\sd_P(f)$ is a weak-equivalence. By \cite[Corollary 1.3.16(b)]{Hovey}, we only need to prove that for any fibrant object in $\sS_P$, $\fil{X}$, the co-unit $\epsilon\colon\sd_P(\Ex_P\fil{X})\to\fil{X}$ is a weak-equivalence. Consider the commutative diagram
\begin{equation*}
\begin{tikzcd}
\sd_P\fil{X}
\arrow{rr}{\lv_P}
\arrow{dr}{\sd_P(\beta)}
&&\fil{X}
\\
&\sd_P(\Ex_P\fil{X})
\arrow{ur}{\epsilon}
\end{tikzcd}
\end{equation*}
As we have seen, the top map is a weak-equivalence. Furthermore, since $\fil{X}$ is fibrant, by \cite[Lemma A.5]{ArticleMoi}, the map $\beta\colon \fil{X}\to\Ex_P\fil{X}$ is a weak-equivalence. Since $\sd_P$ preserve weak-equivalences, $\sd_P(\beta)$ is also a weak-equivalence, and by two out of three, so is $\epsilon\colon \sd_P(\Ex_P\fil{X})\to\fil{X}$.
%
\end{proof}

\section{Vertical filtered CW-complexes and the homotopy theory of filtered spaces}
\label{SectionTopNP}

The goal of this section is to prove that the Quillen adjunction $\varphi_P\circ-\colon \Top_{N(P)}\leftrightarrow \Top_P\colon -\times_P\Real{N(P)}$ is a Quillen equivalence, thereby fixing a gap in the proof of \cite[Theorem 2.15]{TopNP}. 
In \cite{TopNP}, we constructed model categories $\Top_P$ and $\Top_{N(P)}$ related by a Quillen adjunction and showed that it was a Quillen equivalence. The proof of that last claim was incomplete, and so we give a full proof here. By \cite[Corollary 1.3.16(c)]{Hovey}, it suffices to show that the unit of this adjunction, evaluated at a cofibrant object $\fil{X}\in \Top_{N(P)}$ is a weak-equivalence. This is the content of Theorem \ref{TheoremeAMontrer}, which is proven in subsection \ref{SubsectionProvingQETopNP}. In order to prove this theorem, in subsection \ref{SubsectionCW}, we introduce the notion of vertical filtered CW-complexes, (Definition \ref{DefVerticalFilteredCW}). We show that any cofibrant object in $\Top_{N(P)}$ is (strongly) homotopy equivalent to a vertical filtered CW-complex (Proposition \ref{PropCofibrantVerticalCW}), allowing us to reduce the proof of Theorem \ref{TheoremeAMontrer} to such objects.

\subsection{Some results on vertical filtered CW-complexes}
\label{SubsectionCW}
\begin{defin}\label{DefVerticalFilteredCW}
For $n\geq 0$, and $\Delta^{\varphi}\subset N(P)$ a non-degenerate filtered simplex, the filtered cell of dimension $(n,\varphi)$ is the strongly filtered space $(B^n\otimes\Delta^{\varphi},\pr_{\Delta^{\varphi}})$, where $B^n$ is the ball in dimension $n$. Its boundary is the subspace $S^{n-1}\otimes\Delta^{\varphi}$, with the induced filtration.
A filtered CW-complex is a strongly filtered space which can be obtained by inductively gluing cells along their boundaries. More precisely. $\fil{X}$ is a filtered $CW$-complex if a skeletal decomposition of $X$ has been chosen :
\begin{equation*}
X^0\subset X^1\subset \dots\subset X^n\subset\dots \subset \fil{X}
\end{equation*}
such that :
\begin{itemize}
\item $X=\bigcup\limits_{n\in \mathbb{N}}X^n$.
\item $X^0$ is a disjoint union of $0$-cells. (That is, cells of the form $\Delta^{\varphi}$, for $\Delta^{\varphi}\in N(P)$ of any dimension).
\item For all $n\geq 0$, $X^{n+1}$ can be obtained as the following pushout :
\begin{equation*}
\begin{tikzcd}
\coprod\limits_{e_{\alpha}}S^n\otimes\Delta^{\varphi_{\alpha}}
\arrow{d}
\arrow{r}{\coprod \chi_{\alpha}}
&X^n
\arrow{d}
\\
\coprod\limits_{e_{\alpha}}B^{n+1}\otimes\Delta^{\varphi_{\alpha}}
\arrow{r}
&X^{n+1}
\end{tikzcd}
\end{equation*}
where all maps are strongly filtered.
\end{itemize}
A filtered CW-complex is \textbf{vertical} if in addition, for any cell $e_{\alpha}$, its attaching map $\chi_{\alpha}\colon S^{n_{\alpha}}\otimes\Delta^{\varphi_{\alpha}}\to X^{n_{\alpha}}$ satisfy the following condition : If $(x,t)\in S^{n_{\alpha}}\otimes\Delta^{\varphi_{\alpha}}$, and $\chi_{\alpha}(x,t)$ lies in the interior of a cell $e_{\beta}$ of dimension $(m,\Delta^{\psi})$, and has coordinates $((\chi_{\alpha})_1(x,t),(\chi_{\alpha})_2(x,t))$ in that cell. then :
\begin{itemize}
\item $\Delta^{\varphi_{\alpha}}\subset\Delta^{\psi}$
\item For all $s\in \Delta^{\varphi_{\alpha}}$, $\chi_{\alpha}(x,s)$ lies in the interior of the cell $e_{\beta}$ and has coordinates $((\chi_{\alpha})_1(x,t),s)$. 
\end{itemize}

A filtered map $f\colon \fil{X}\to\fil{Y}$ between vertical filtered CW-complexes is vertical if for any cell of $X$, $e_{\alpha}$, and any cell of $Y$, $e_{\beta}$ whose interior intersects $f(e_{\alpha})$, the following conditions are satisfied :
\begin{itemize}
\item $\Delta^{\varphi_{\alpha}}\subset \Delta^{\varphi_{\beta}}$
\item For $(x,t)\in B^{n_{\alpha}}\otimes \Delta^{\varphi_{\alpha}}$, if $f_{\alpha}(x,t)=(y,t)\in e_{\beta}$, then for all $s\in \Delta^{\varphi_{\alpha}}$, $f_{\alpha}(x,s)=(y,s)$.
\end{itemize}
We write $\Ver_P$ for the category of vertical filtered CW-complexes over $P$ and vertical filtered maps.
\end{defin}

\begin{construction}\label{RemarqueConstructionCW}
To any vertical filtered CW-complex, $\fil{X}$, we associate a (non-filtered) CW-complex, $L(X)$, in the following way.
Consider some attaching map $\chi_{\alpha}\colon S^{n_{\alpha}}\otimes \Delta^{\varphi_{\alpha}}\to X^{n_{\alpha}}$. Pick a cell, $e_{\beta}$ of $X$ whose interior intersects the image of $\chi_{\alpha}$, and consider the restriction
\begin{equation*}
\chi_{\alpha}\colon\chi_{\alpha}^{-1}(e_{\beta})\to B^{n_{\beta}}\otimes \Delta^{\varphi_{\beta}}.
\end{equation*}
The verticality hypothesis guarantees that on this restriction, $\chi_{\alpha}$ can be decomposed as the product of two maps :
\begin{itemize}
\item a continuous map $L(\chi_{\alpha}^{\beta})$ from a subset of $\colon S^{n_{\alpha}}$ to $B^{n_{\beta}}$
\item the strongly filtered inclusion $\Delta^{\varphi_{\alpha}}\hookrightarrow \Delta^{\varphi_{\beta}}$
\end{itemize}
We construct a non-filtered CW-complex, $L(X)$, by giving it a $n$-cell, $L(e_{\alpha})$, for each $(n,\Delta^{\varphi})$-cell, $e_{\alpha}$, of $X$, and building the attaching maps from the $L(\chi_{\alpha}^{\beta})$. More explicitly, we define $L(X)$ by induction. Set 
\begin{equation*}
L(X)^0=\coprod\limits_{e_{\alpha}, \dim(e_{\alpha})=(0,\Delta^{\varphi_{\alpha}})}B^0
\end{equation*}
Then, assume that $L(X)^n$ has been built, and let $e_{\alpha}$ be a cell of dimension $(n+1,\Delta^{\varphi_{\alpha}})$. 
As we have seen there are well defined maps $L(\chi_{\alpha}^{\beta})$ from subsets of $S^n$ to the cells, $L(e_{\beta})$, of $L(X)^n$. Since $\chi_{\alpha}$ is an attaching map, the pre-images $\chi_{\alpha}^{-1}(e_{\beta})$, for all cells $e_{\beta}\in X^n$, must cover $S^n\otimes \Delta^{\varphi_{\alpha}}$. Since $\fil{X}$ is a vertical filtered CW-complex, this implies that the codomains of the $L(\chi_{\alpha}^{\beta})$, must cover $S^n$, and that they coincide on the intersection of the codomains. In turn, they can be glued to produce an attaching map for $L(e_{\alpha})$, $L(\chi_{\alpha})\colon S^n\to  L(X)^n$.
We can now define $L(X)^{n+1}$ as the following pushout :
\begin{equation*}
\begin{tikzcd}
\coprod\limits_{e_{\alpha},\dim(e_{\alpha})=(n+1,\Delta^{\varphi_{\alpha}})}S^n
\arrow{r}{\coprod L(\chi_{\alpha})}
\arrow{d}
&L(X)^n
\arrow{d}
\\
\coprod\limits_{e_{\alpha},\dim(e_{\alpha})=(n+1,\Delta^{\varphi_{\alpha}})}B^{n+1}
\arrow{r}
&L(X)^{n+1}
\end{tikzcd}
\end{equation*}
 Furthermore, assume $f\colon \fil{X}\to\fil{Y}$ is a vertical filtered map between vertical filtered CW-complexes. Just as earlier, if $e_{\alpha}$ is a cell of $X$, and $e_{\beta}$ is a cell of $Y$ whose interior intersects $f(e_{\alpha})$, the verticality condition allows us to decompose a suitable restriction of $f$ as the product of two maps :
\begin{itemize}
\item a continuous map $L(f_{\alpha}^{\beta})$ from a subset of $B^{n_{\alpha}}$ to $B^{n_{\beta}}$,
\item the strongly filtered inclusion $\Delta^{\varphi_{\alpha}}\hookrightarrow \Delta^{\varphi_{\beta}}$.
\end{itemize}
First, one can glue the $L(f_{\alpha}^{\beta})$ for all the cells $e_{\beta}$ intersecting $f(e_{\alpha})$. This gives a map $L(f_{\alpha})\colon B^{n_{\alpha}}\to L(Y)$. One can then glue those $L(f_{\alpha})$ together, to get a map $L(f)\colon L(X)\to L(Y)$.
\end{construction}

We have just defined a functor from the category of vertical filtered CW-complexes to the category of CW-complexes. We will see that, provided we keep the information about the filtered dimensions of the cells, this functor retains all the information about vertical filtered CW-complexes. To make this statement precise, we need the following definition.

\begin{defin}
Let $K$ be a CW-complex. A $P$-labeling of $K$ is the data of a map $\lambda_K$ from the set of cells of $K$ to the set of non-degenerate simplices of $N(P)$ such that, if $e_{\alpha}$ intersects the closure of $e_{\beta}$, then $\lambda_K(e_{\beta})\subset \lambda_{K}(e_{\alpha})$. A label-preserving map between $P$-labeled CW-complexes $f\colon \lab{K}\to\lab{L}$ is a continuous map $f\colon K\to L$ such that for any cell of $K$, $e_{\alpha}$ and any cell of $L$, $e_{\beta}$, if $f(e_{\alpha})$ intersects $e_{\beta}$, then $\lambda_K(e_{\alpha})\subset \lambda_L(e_{\beta})$. Let $\CW_P$ be the category of $P$-labeled CW-complexes and label-preserving maps.
\end{defin}

\begin{prop}\label{PropositionEquivalenceVertCWP}
The categories $\Ver_P$ and $CW_P$ are equivalent.
\end{prop}

\begin{proof}
We will construct both functors in this equivalence.
In construction \ref{RemarqueConstructionCW}, we constructed a functor $L\colon\Ver_P\to \CW$. We need to define a $P$-labeling on its image. Let $\fil{X}$ be a vertical filtered CW-complex. The CW-complex $L(X)$ contains a cell $L(e_{\alpha})$ for each cell $e_{\alpha}$ of $X$. 
If $e_{\alpha}$ is a cell of dimension $(n_{\alpha},\Delta^{\varphi_{\alpha}})$, we define the labeling on $L(e_{\alpha})$ as $\lambda_{L(X)}(L(e_{\alpha}))=\Delta^{\varphi_{\alpha}}$. 
We need to check that the map $\lambda_{L(X)}$ satisfies the condition of being a $P$-labeling. 
Assume that $L(e_{\alpha})$ intersects the closure of $L(e_{\beta})$. 
This means that $L(e_{\alpha})$ intersects the image of the attaching map 
$L(\chi_{\beta})\colon S^{n_{\beta}-1}\to L(X)^{n_{\beta}-1}$. 
By construction of $L(X)$, this implies that $e_{\alpha}$ intersects the image of the attaching map of $e_{\beta}$, $\chi_{\beta}\colon S^{n_{\beta}-1}\otimes\Delta^{\varphi_{\beta}}\to X^{n_{\beta}-1}$. But, by the verticality condition, this implies that $\Delta^{\varphi_{\beta}}\subset \Delta^{\varphi_{\alpha}}$.
The same check at the level of maps shows that the association $\fil{X}\to\lab{L(X)}$ induces a well-defined functor $\Ver_P\to \CW_P$. By abuse of notation, we will still write $L$ for this functor. From now on, $L\fil{X}$ denotes the $P$-labeled CW-complex, while $L(X)$ denotes the underlying CW-complex (without labels).

We still need to construct the inverse functor. Let $\lab{K}$ be a $P$-labeled CW-complex. We will define a vertical filtered CW-complex $\fil{V(K)}$ with one $(n,\Delta^{\varphi})$-cell for every $n$-cell of $K$ with label $\Delta^{\varphi}$. Define the $0$-skeleton of $V(K)$ as 
\begin{equation*}
V(K)^0=\coprod\limits_{e_{\alpha}, \dim(e_{\alpha})=0} \lambda_K(e_{\alpha}).
\end{equation*}
Then, assume that $V(K)^{n}$ has been built, and let $e_{\alpha}$ be a $(n+1)$-cell of $K$. The attaching map $\chi_{\alpha}\colon S^n\to K^n$ lands in cells $e_{\beta}$ all satisfying $\lambda_K(e_{\alpha})\subset \lambda_K(e_{\beta})$. 
In addition, for any such $e_{\beta}$, there is a well defined restriction $\chi_{\alpha}^{\beta}\colon \chi_{\alpha}^{-1}(e_{\beta})\to e_{\beta}$. 
Since the cell $V(e_{\beta})$ is in $V(K)^n$, and of dimension $(n_{\beta},\lambda_K(e_{\beta}))$, we can define a map $V(\chi_{\alpha}^{\beta})\colon \chi_{\alpha}^{-1}(e_{\beta})\otimes\lambda_K(e_{\alpha})\to V(e_{\beta})$ as the product of the map 
$\chi_{\alpha}^{\beta}$, and the inclusion 
$\lambda_K(e_{\alpha})\hookrightarrow \lambda_K(e_{\beta})$. 
Gluing the $V(\chi_{\alpha}^{\beta})$ together, we get a vertical filtered map 
$V(\chi_{\alpha})\colon S^n\otimes \lambda_K(e_{\alpha})\to V(K)^n$. 
We can now define $V(K)^{n+1}$ as the following push-out :
\begin{equation*}
\begin{tikzcd}
\coprod\limits_{e_{\alpha}, \dim(e_{\alpha})=n+1} S^{n}\otimes\lambda_K(e_{\alpha})
\arrow{d}
\arrow{r}{\coprod V(\chi_{\alpha})}
&V(K)^n
\arrow{d}
\\
\coprod\limits_{e_{\alpha}, \dim(e_{\alpha})=n+1} B^{n+1}\otimes\lambda_K(e_{\alpha})
\arrow{r}
&V(K)^{n+1}
\end{tikzcd}
\end{equation*}
The proof that $V$ sends label-preserving maps to vertical maps is similar. This imply that $V\colon \CW_P\to \Ver_P$ is a well-defined functor. But by construction, $V(L\fil{X})\simeq \fil{X}$ and $L(V\lab{K})\simeq \lab{K}$, which concludes the proof.
\end{proof}

\begin{remarque}
The equivalence of categories of proposition \ref{PropositionEquivalenceVertCWP} is compatible with homotopies in the following sense : two vertical maps $f,g\colon \fil{X}\to\fil{Y}$ are homotopic through a vertical homotopy if and only if $V(f)$ and $V(g)$ are homotopic through a homotopy preserving labels. This follows from the fact that $V$ and $L$ preserve cylinders.
\end{remarque}

\begin{remarque}
Let $\lab{K}$ be a $P$-labeled CW-complex, and $U\subset K$ a subspace. One can define a subspace $V(U,\lambda_K)\subset V\lab{K}$ in the following way. For any cell $e_{\alpha}\in K$, consider the subspace $(e_{\alpha}\cap U)\otimes\lambda_K(e_{\alpha})\subset V(e_{\alpha})$. The union of those subspaces is the subspace $V(U,\lambda_K)$. Furthermore, if we assume $U\subset K$ to be open, then $V(U,\lambda_K)\subset V\lab{K}$ is an open subset.
\end{remarque}

\begin{defin}
A vertical map between vertical filtered CW-complexes $f\colon\fil{X}\to\fil{Y}$ is said to be cellular if for all cell of $X$, $e_{\alpha}$, $f(e_{\alpha})\subset Y^{n_{\alpha}}$.
\end{defin}

We will make use of the following few lemmas

\begin{lemme}\label{LemmePushoutSingleCell}
Let $\fil{X}$ be a vertical filtered CW-complex, and $\chi\colon S^n\otimes\Delta^{\varphi}\to \fil{X}$ be a cellular vertical map. Define $\fil{Y}=\fil{X}\cup_{\chi} (B^{n+1}\otimes \Delta^{\varphi})$ as the following pushout.
\begin{equation*}
\begin{tikzcd}
S^n\otimes\Delta^{\varphi}
\arrow{r}{\chi}
\arrow{d}
&\fil{X}
\arrow{d}
\\
B^{n+1}\otimes\Delta^{\varphi}
\arrow{r}
&\fil{Y}
\end{tikzcd}
\end{equation*}
Then, $\fil{Y}$ is a vertical filtered CW-complex, and $\fil{X}$ is a subcomplex of $\fil{Y}$.
\end{lemme}

\begin{proof}
By assumption, $\fil{X}$ is a vertical filtered CW-complex, which means it admits a skeletal decomposition :
\begin{equation*}
X^0\subset X^1\subset \dots \subset X^n\subset \dots \subset X
\end{equation*}
We define the following skeletal decomposition for $Y$.
For $k\leq n$, let $Y^k=X^k$. Then, for $k>n$, let $Y^k=X^k\cup_{\chi} B^{n+1}\otimes\Delta^{\varphi}$. Since $\chi$ is supposed to be cellular, we have $\Im(\chi)\subset X^n$, and the $Y^k$ are all well-defined, and contain $X^k$. Furthermore, for $k\not= n$, $Y^{k+1}$ is obtained as the following pushout
\begin{equation*}
\begin{tikzcd}
\coprod\limits_{\substack{e_{\alpha}\in X \\ \dim(e_{\alpha})=(k+1,\Delta^{\varphi_{\alpha}})}}S^k\otimes\Delta^{\varphi_{\alpha}}
\arrow{r}{\coprod\chi_{\alpha}}
\arrow{d}
&Y^k
\arrow{d}
\\
\coprod\limits_{\substack{e_{\alpha}\in X \\ \dim(e_{\alpha})=(k+1,\Delta^{\varphi_{\alpha}})}}B^{k+1}\otimes\Delta^{\varphi_{\alpha}}
\arrow{r}
&Y^{k+1}
\end{tikzcd}
\end{equation*}
And, for $k=n$, $Y^{n+1}$ can be obtained as the following pushout :
\begin{equation*}
\begin{tikzcd}
S^n\otimes\Delta^{\varphi}\coprod \coprod\limits_{\substack{e_{\alpha}\in X \\ \dim(e_{\alpha})=(n+1,\Delta^{\varphi_{\alpha}})}}S^n\otimes\Delta^{\varphi_{\alpha}}
\arrow{r}{\chi\coprod \coprod\limits_{e_{\alpha}}\chi_{\alpha}}
\arrow{d}
&Y^n
\arrow{d}
\\
B^{n+1}\otimes\Delta^{\varphi}\coprod \coprod\limits_{\substack{e_{\alpha}\in X \\ \dim(e_{\alpha})=(n+1,\Delta^{\varphi_{\alpha}})}}B^{n+1}\otimes\Delta^{\varphi_{\alpha}}
\arrow{r}
&Y^{n+1}
\end{tikzcd}
\end{equation*}
This proves that $\fil{Y}$ is a vertical filtered CW-complex. It clearly contains $\fil{X}$ as a subspace, and by construction, cells of $\fil{X}$ are also cells of $\fil{Y}$, which means that $\fil{X}$ is a subcomplex.
\end{proof}

\begin{lemme}\label{LemmeCWPushout}
Let $\fil{A},\fil{X}$ and $\fil{Z}$ be three vertical filtered CW-complexes, and let $f\colon \fil{A}\to \fil{Z}$ be a cellular vertical map, and $i\colon \fil{A}\to \fil{X}$ the inclusion of a subcomplex (that is, all cells of $A$ are cells of $X$). Assume $\fil{W}$ is the following pushout
\begin{equation*}
\begin{tikzcd}
\fil{A}
\arrow{r}{f}
\arrow[swap]{d}{i}
&\fil{Z}
\arrow{d}{j}
\\
\fil{X}
\arrow{r}{g}
&\fil{W}
\end{tikzcd}
\end{equation*}
Then, $\fil{W}$ admits a structure of vertical filtered CW-complex for which $j$ is the inclusion of a sub-complex, and $g$ is a cellular vertical map.
\end{lemme}

\begin{proof}
Consider the filtration of $X$, $A\cup X^0\subset A\cup X^1\dots \subset X$. For all $n\geq 0$, we have the push-out
\begin{equation*}
\begin{tikzcd}
\coprod\limits_{\substack{e_{\alpha}\in X\setminus A\\ \dim(e_{\alpha})=(n+1,\Delta^{\varphi_{\alpha}})}}S^n\otimes\Delta^{\varphi_{\alpha}}
\arrow{r}{\coprod\chi_{\alpha}}
\arrow{d}
&A\cup X^n
\arrow{d}
\\
\coprod\limits_{\substack{e_{\alpha}\in X\setminus A\\ \dim(e_{\alpha})=(n+1,\Delta^{\varphi_{\alpha}})}}B^{n+1}\otimes\Delta^{\varphi_{\alpha}}
\arrow{r}
&A\cup X^{n+1}
\end{tikzcd}
\end{equation*}
Now define $Z\cup_A X^0$ as $Z\cup\coprod\limits_{\substack{e_{\alpha}\in X\setminus A\\ \dim(e_{\alpha})=(0,\Delta^{\varphi_{\alpha}})}}\Delta^{\varphi_{\alpha}}$. By construction, there exists a map $f^0\colon A\cup X^0\to Z\cup_A X^0$ extending $f$. Now assume that $Z\cup_A X^n$ has been constructed, together with a cellular vertical map $f^n\colon A\cup X^n\to Z\cup_A X^n$ and define $Z\cup_A X^{n+1}$ as the following push-out
\begin{equation*}
\begin{tikzcd}
\coprod\limits_{\substack{e_{\alpha}\in X\setminus A\\ \dim(e_{\alpha})=(n+1,\Delta^{\varphi_{\alpha}})}}S^n\otimes\Delta^{\varphi_{\alpha}}
\arrow{r}{\coprod\chi_{\alpha}}
\arrow{d}
&A\cup X^n
\arrow{d}
\arrow{r}{f^n}
&Z\cup_A X^n
\arrow{d}
\\
\coprod\limits_{\substack{e_{\alpha}\in X\setminus A\\ \dim(e_{\alpha})=(n+1,\Delta^{\varphi_{\alpha}})}}B^{n+1}\otimes\Delta^{\varphi_{\alpha}}
\arrow{r}
&A\cup X^{n+1}
\arrow{r}{f^{n+1}}
&Z\cup_A X^{n+1}
\end{tikzcd}
\end{equation*}
Since, $\fil{X}$ is a filtered vertical CW-complex, the attaching maps $\chi_{\alpha}$ must land in the $n$-skeleton of $X$. But since $f^n$ is cellular, this implies that for any cell of dimension $n+1$, the attaching maps $f\circ \chi_{\alpha}$ must land in the $n$-skeleton of $Z\cup_A X^n$. In particular, by lemma \ref{LemmePushoutSingleCell}, $Z\cup_A X^{n+1}$ is a filtered vertical CW-complex. Furthermore, since by assumption $f^n$ was vertical and cellular, and since $f^{n+1}$ is the identity when restricted to the interior of the newly added cells, the map $f^{n+1}$ is still vertical and cellular. We then conclude by lemma \ref{LemmeColimiteCW}.
\end{proof}

\begin{lemme}\label{LemmeColimiteCW}
Let $(Y_n)_{n\in \N}$ be a family of vertical filtered CW-complexes together with maps $i_n\colon Y_n\hookrightarrow Y_{n+1}$ that are inclusions of sub-complexes. Then $Y=\colim Y_n$ is a vertical filtered CW-complex, and the maps $Y_n\to Y$ are inclusions of sub-complexes for all $n\geq 0$.
\end{lemme}

\begin{proof}
The filtered space $Y$ admits a vertical CW-structure where it has a cell for each cell eventually appearing in the sequence $Y_0\subset Y_1\subset ...$. This structure turns all the $Y_n$ into sub-complexes.
\end{proof}

\begin{prop}\label{PropRealSdCW}
Let $\fil{X}$ be a filtered simplicial set. The realisation of its subdivision, $\RealNP{\sd_P\fil{X}}$ admits the structure of a vertical filtered CW-complex.
\end{prop}

\begin{proof}
Given a filtered simplicial set $X$, One can consider the filtration by sub-simplicial sets $X^n\subset X$ generated by simplices of dimension $\leq n$. This gives rise to the description of $\RealNP{\sd_P\fil{X}}$ as the colimit $\colim_n\RealNP{\sd_P\fil{X^n}}$. By lemma \ref{LemmeColimiteCW}, it suffices to show that $\RealNP{\sd_P(X^0,\varphi_{X^0})}$ is a vertical filtered CW-complex (this is clear, since it is a disjoint union of points), and that for all $n\geq 0$, the map $\RealNP{\sd_P\fil{X^n}}\to \RealNP{\sd_P\fil{X^{n+1}}}$ is the inclusion of a sub-complex.  Note that for all $n\geq 0$, $\RealNP{\sd_P\fil{X^{n+1}}}$ can be described as the following push-out :
\begin{equation*}
\begin{tikzcd}[column sep = 80]
\coprod\limits_{\sigma\in X_{n+1}^{\nd}}\RealNP{\sd_P(\partial\Delta^{\varphi_{\sigma}})}
\arrow{r}{\coprod (\RealNP{\sd_P(\sigma)})_{|\partial\Delta^{\varphi_{\sigma}}}}
\arrow{d}
&\RealNP{\sd_P\fil{X^n}}
\arrow{d}
\\
\coprod\limits_{\sigma\in X_{n+1}^{\nd}}\RealNP{\sd_P(\Delta^{\varphi_{\sigma}})}
\arrow{r}
&\RealNP{\sd_P\fil{X^{n+1}}}
\end{tikzcd}
\end{equation*}
Working by induction and using lemma \ref{LemmeCWPushout}, it is enough to show that the top map is cellular and vertical, and that the leftmost map is the inclusion of a sub-complex. Since those properties are preserved by taking disjoint unions, we have to show that :
\begin{itemize}
\item For any filtered simplex, $\Delta^{\varphi}$, $\RealNP{\sd_P(\Delta^{\varphi})}$ is a vertical filtered CW-complex, and $\RealNP{\sd_P(\partial\Delta^{\varphi})}\to \RealNP{\sd_P(\Delta^{\varphi})}$ is the inclusion of a sub-complex. 
\item for any filtered simplex of $X$, $\sigma\colon \Delta^{\varphi_{\sigma}}\to X$, of dimension $n+1$, the map
\begin{equation}\label{EquationSdSigmaDansX}
\RealNP{\sd_P(\partial(\Delta^{\varphi_{\sigma}}))}\to \RealNP{\sd_P\fil{X^n}}
\end{equation}
is cellular and vertical.
\end{itemize}
The former is covered in the first part of Lemma \ref{LemmeSdSimplexVerticalCW}. For the latter, consider a filtered simplex of $X$, $\sigma\colon \Delta^{\varphi_{\sigma}}\to X$ of dimension $n+1$, and consider its restriction to a face $d_i\Delta^{\varphi_{\sigma}}\to X$ for some $0\leq i\leq n+1$. By definition, this restriction is a filtered simplex of dimension $n$ of $X$, $\tau\colon \Delta^{\varphi_{\tau}}\to X$. In particular, it realizes to a map
\begin{equation*}
\RealNP{\sd_P(d_i\Delta^{\varphi_{\sigma}})}=\RealNP{\sd_P(\Delta^{\varphi_{\tau}})}\to \RealNP{\sd_P\fil{X^n}}
\end{equation*}
Which can be assumed to be cellular and vertical, by the induction hypothesis. Since cellularity and verticality are local conditions, it suffices to check that they are satisified for all cells of $\RealNP{\sd_P(\partial(\Delta^{\varphi_{\sigma}}))}$ to conclude that \eqref{EquationSdSigmaDansX} is cellular and vertical. But all those cells appear in at least one of the $\RealNP{\sd_P(d_i\Delta^{\varphi_{\sigma}})}$. We deduce that the map \eqref{EquationSdSigmaDansX} is vertical and cellular. Furthermore, the map induced by $\sigma$ :
\begin{equation*}
\RealNP{\sd_P(\Delta^{\varphi_{\sigma}})}\to \RealNP{\sd_P\fil{X^{n+1}}}
\end{equation*}
is the push-out of the vertical and cellular map \eqref{EquationSdSigmaDansX} along a sub-complex inclusion. By lemma \ref{LemmeCWPushout}, this implies that it is a vertical and cellular map, concluding the proof.
\end{proof}

\begin{lemme}\label{LemmeSdSimplexVerticalCW}
Let $\Delta^{\varphi}$ be a filtered simplex. $\RealNP{\sd_P(\Delta^{\varphi})}$ admits a structure of vertical filtered CW-complex, such that $\RealNP{\sd_P(\partial\Delta^{\varphi})}$ and $\RealNP{\sd_P(d_i\Delta^{\varphi})}$ are subcomplexes.
%
\end{lemme}

\begin{proof}
Throughout this proof, we will identify $S^n$ with $\Real{\partial\Delta^{n+1}}$ and $B^{n+1}$ with $\Real{\Delta^{n+1}}$.
Let $\Delta^{\varphi}$ be a filtered simplex. We will define a vertical filtered CW-complex $X$, and show that it is filtered homeomorphic to $\RealNP{\sd_P(\Delta^{\varphi})}$. Recall that if $\Delta^{\psi}$ is a possibly degenerate filtered simplex, $\Delta^{\overline{\psi}}$ is the unique non-degenerate simplex of which $\Delta^{\psi}$ is a degeneracy. $X$ will have a cell of dimension $(k,\Delta^{\overline{\psi_0}})$, $e_{\psi_0,\dots,\psi_k}$, for each strictly increasing chain $\Delta^{\psi_0}\subsetneq \Delta^{\psi_1}\subsetneq\dots\subsetneq \Delta^{\psi_k}\subset \Delta^{\varphi}$.  Define $X^0$ as
\begin{equation*}
X^0=\coprod\limits_{\Delta^{\psi}\subset\Delta^{\varphi}} \Delta^{\overline{\psi}}
\end{equation*}
Now, let $n\geq 0$, and assume that $X^n$ has been built, with a $(k,\Delta^{\overline{\psi_0}})$-cell for every chain $\Delta^{\psi_0}\subsetneq\Delta^{\psi_1}\subsetneq\dots\subsetneq \Delta^{\psi_k}\subset \Delta^{\varphi}$, for $k\leq n$. Consider some chain $\Delta^{\psi_0}\subsetneq\Delta^{\psi_1}\subsetneq\dots\subsetneq \Delta^{\psi_{k+1}}\subset \Delta^{\varphi}$. We need an attaching map $\chi_{\psi_0,\dots,\psi_{k+1}}\colon S^n\otimes \Delta^{\overline{\psi_0}}\to X^n$ to glue the corresponding cell. Identify $S^n$ with $\Real{\partial\Delta^{n+1}}$.
We can then cut the boundary of the cell $e_{\psi_0,\dots,\psi_{k+1}}$ into the subspaces $\Real{d_i\Delta^{n+1}}\otimes \Delta^{\overline{\psi_0}}$, for $0\leq i\leq n+1$.
For $i\not = 0$, the cells $e_{\psi_0,\dots,\widehat{\psi_i},\dots,\psi_{k+1}}$ are of dimension $(n,\Delta^{\overline{\psi_0}})$.
In particular, consider the homeomorphism $(D^i)^{-1}\colon \Real{d_i\Delta^{n+1}}\to \Real{\Delta^n}$. Taking the product of this homeomorphism, with the identity $\Id\colon \Delta^{\overline{\psi_0}}\to\Delta^{\overline{\psi_0}}$ gives part of the attaching map 
\begin{equation*}
\chi_{\psi_0,\dots,\psi_{k+1}}^i\colon\Real{d_i\Delta^{n+1}}\otimes \Delta^{\overline{\psi_0}}\simeq \Real{\Delta^n}\otimes\Delta^{\overline{\psi_0}}\xrightarrow{e_{\psi_0,\dots,\widehat{\psi_i},\dots,\psi_n}} X^n
\end{equation*}
In the case $i=0$, the cell $e_{\psi_1,\dots,\psi_{n+1}}$ is of dimension $(n,\overline{\psi_1})$, with $\Delta^{\overline{\psi_0}}\subseteq \Delta^{\overline{\psi_1}}$. 
This means that we can  similarly define a vertical map from $\Real{d_0\Delta^{n+1}}\otimes \Delta^{\overline{\psi_0}}$ 
to the closure of $e_{\psi_1,\dots,\psi_{n+1}}$, by taking the product of 
$(D^{0})^{-1}$ with the inclusion 
$\Delta^{\overline{\psi_0}}\hookrightarrow\Delta^{\overline{\psi_1}}$. 
By construction, on any intersection $d_i\Delta^{n+1}\cap d_j\Delta^{n+1}$, the partial attaching maps $\chi^i$ and $\chi^j$ coincide, meaning they can be glued to yield an attaching map $\chi_{\psi_0,\dots,\psi_{n+1}}\colon \Real{\partial\Delta^{n+1}}\otimes\Delta^{\overline{\psi_0}}\to X^n$.
We can then define $X^{n+1}$ as the following push-out:
\begin{equation*}
\begin{tikzcd}[column sep = 60]
\coprod\limits_{\Delta^{\psi_0}\subsetneq\dots\subsetneq\Delta^{\psi_{n+1}}\subset \Delta^{\varphi}} S^n\otimes \Delta^{\overline{\psi_0}}
\arrow{r}{\coprod \chi_{\psi_0,\dots,\psi_{n+1}}}
\arrow{d}
&X^n
\arrow{d}
\\
\coprod\limits_{\Delta^{\psi_0}\subsetneq\dots\subsetneq\Delta^{\psi_{n+1}}\subset \Delta^{\varphi}} B^{n+1}\otimes \Delta^{\overline{\psi_0}}
\arrow{r}
&X^{n+1}
\end{tikzcd}
\end{equation*}
It remains to be shown that $\fil{X}$ is filtered homeomorphic to $\RealNP{\sd_P(\Delta^{\varphi})}$. Recall that $\sd_P(\Delta^{\varphi})$ can be defined as the following pullback
\begin{equation*}
\begin{tikzcd}
\sd_P(\Delta^{\varphi})
\arrow{r}
\arrow[swap]{d}{\sd_P(\varphi)}
&\sd(\Delta^{\varphi})
\arrow{d}{\sd(\varphi)}
\\
\sd_P(N(P))
\arrow{r}
&\sd(N(P))
\end{tikzcd}
\end{equation*}
Since the functor $\RealNP{-}$ preserves pullback, it suffices to show that $X$ fits in a pullback square
\begin{equation}\label{EquationPullbackCWsdP}
\begin{tikzcd}
X
\arrow{r}{f}
\arrow[swap]{d}{g}
&\Real{\sd(\Delta^{\varphi})}
\arrow{d}
\\
\Real{\sd_P(N(P))}
\arrow{r}
&\Real{\sd(N(P))}
\end{tikzcd}
\end{equation}
and that the stratification on $X$ is given by the composition $X\to \Real{\sd_P(N(P))}\to \Real{N(P)}$. First, note that by construction, $\Real{\sd(\Delta^{\varphi})}$ is homeomorphic to $L(X)$ (as defined in remark \ref{RemarqueConstructionCW}). In particular, $f$ will send the cells $e_{\psi_0,\dots,\psi_n}$ to the corresponding cells $L(e_{\psi_0,\dots,\psi_n})$. And, for any point in this cell, $(x,t)\in B^n\otimes\Delta^{\overline{\psi_0}}$, set $f_{\psi_0,\dots,\psi_k}(x,t)=x\in L(e_{\psi_0,\dots,\psi_n})$. 

Now to define $g$, recall that $\sd_P(N(P))$ is a simplicial complex with a vertex for each pair $(\Delta^{\psi},p)$, where $\Delta^{\psi}$ is a non-degenerate filtered simplex and $p\in \Delta^{\psi}$. The higher dimensional simplices of $\sd_P(N(P))$ are then given by the tuples $\{(\Delta^{\psi_0},p_0),\dots,(\Delta^{\psi_n},p_n)\}$, where $\Delta^{\psi_0}\subset\dots\subset \Delta^{\psi_n}$ is an increasing chain of non-degenerate simplices, and $p_i\in \Delta^{\psi_0}$, for all $0\leq i\leq n$. Consider a cell in $\fil{X}$, $e_{\psi_0,\dots,\psi_n}$. Topologically, it is a product of two simplices, $\Delta^n$ and $\Delta^{\overline{\psi_0}}$. As such, one can triangulate it in a canonical way. Label the vertices of $\Delta^n$ by the $\Delta^{\psi_i}$ and write $[p_0,\dots,p_k]=\Delta^{\overline{\psi_0}}$. The vertices of the triangulation of $\Delta^n\otimes \Delta^{\overline{\psi_0}}$ will then be given by the pairs $(\Delta^{\psi_i},p_j)$, for $0\leq i\leq n$, $0\leq j\leq k$. Furthermore, simplices of this triangulation will be given by chains $((\Delta^{\psi_{i_0}},p_{j_0}),\dots,(\Delta^{\psi_{i_l}},p_{j_l}))$, where $\Delta^{\psi_{i_{m}}}\subset \Delta^{\psi_{i_{m+1}}}$, and $p_{j_m}\leq p_{j_{m+1}}$ for all $0\leq m\leq {l-1}$. The isomorphism between this simplicial set and $e_{\psi_0,\dots,\psi_n}$ is then given by sending the vertices of the simplicial sets to the corresponding points in the cell, and then extending linearly. One checks that this also turns the map $f$ defined earlier into a simplicial map.
On the other hand, we can now define $g_{\psi_0,\dots,\psi_n}\colon \Real{\Delta^n}\otimes\Delta^{\overline{\psi_0}}\to \Real{\sd_P(N(P))}$ by setting $g(\Delta^{\psi_i},p_j)=(\Delta^{\overline{\psi_i}},p_j)$ and extending linearly.

It remains to be shown that the square \eqref{EquationPullbackCWsdP} is a pullback square. Let us first prove that it is commutative. Let $e_{\psi_0,\dots,\psi_n}\in X$ be a cell.
 With the previous triangulation of the cell the restriction of $f$ and $g$ to $e_{\psi_0,\dots,\psi_n}$ are simplicial.
  In particular, it is enough to check the commutativity at the level of vertices. Given $(\Delta^{\psi_i},p_j)$ a vertex of the triangulation, we have $f(\Delta^{\psi_i},p_j)=\Delta^{\psi_i}$,
   and its image in $\sd(N(P))$ is $\Delta^{\overline{\psi_i}}$. 
   On the other hand, we have $g(\Delta^{\psi_i},p_j)=(\Delta^{\overline{\psi_i}},p_j)$,
    whose image in $\sd(N(P))$ is $\Delta^{\overline{\psi_i}}$. We conclude that, \eqref{EquationPullbackCWsdP} is commutative. Now, consider a simplex in $\sd(N(P))$, $(\Delta^{\mu_0},\dots,\Delta^{\mu_n})$, with $\Delta^{\mu_0}\subsetneq\dots\subsetneq\Delta^{\mu_n}\subset N(P)$. Its preimage in $\sd_P(N(P))$ is isomorphic as a simplicial complex to $\Delta^n\otimes \Delta^{\mu_0}$. On the other hand, its preimage in $\sd(\Delta^{\varphi})$ is the full sub-simplicial complex spanned by the vertices $\Delta^{\psi}$ such that $\Delta^{\overline{\psi}}=\Delta^{\mu_i}$ for some $0\leq i\leq n$. Finally, its preimage in $X$ is the subspace spanned by all the cells of the form $e_{\psi_0,\dots,\psi_k}$, where for all $0\leq i\leq k$, $\Delta^{\overline{\psi_i}}=\Delta^{\mu_j}$ for some $0\leq j\leq n$. 
    But one notes that such a sell $e_{\psi_0,\dots,\psi_k}\simeq \Delta^k\otimes \Delta^{\overline{\psi_0}}$ is precisely the product of $(\Delta^{\overline{\psi_0}},\dots,\Delta^{\overline{\psi_k}})\otimes\Delta^{\overline{\psi_0}}\subset \sd_P(N(P))$ with $(\Delta^{\psi_0},\dots,\Delta^{\psi_k})\subset \sd(\Delta^{\varphi})$, over $(\Delta^{\mu_0},\dots,\Delta^{\mu_n})$. In particular, the square \eqref{EquationPullbackCWsdP} is indeed a pull-back square. Furthermore, the stratification on $\fil{X}$, which is given on cells of the form $\Real{\Delta^n}\otimes\Delta^{\mu}$ as the projection on $\Delta^{\mu}$, clearly factors through $g$. In particular, $\fil{X}$ is a vertical filtered CW structure on $\RealNP{\sd_P(\Delta^{\varphi})}$. Furthermore, the subspace $\RealNP{\sd_P(\partial\Delta^{\varphi})}$ is spanned by the cells of the form $e_{\psi_0,\dots,\psi_n}$ where $\Delta^{\psi_n}\not = \Delta^{\varphi}$. In particular, it is indeed a sub-complex. Similarily, $\RealNP{\sd_P(d_i\Delta^{\varphi})}$ is the sub-complex containing all cells $e_{\psi_0,\dots,\psi_n}$ such that $\Delta^{\psi_n}\subset d_i\Delta^{\varphi}$.
%
\end{proof}

\begin{prop}\label{PropCofibrantVerticalCW}
Let $\fil{X}$ be a cofibrant object in $\Top_{N(P)}$. There exists a vertical filtered CW-complex, $\fil{Y}$, such that $\fil{X}$ and $\fil{Y}$ are filtered homotopy equivalent in $\Top_{N(P)}$.
\end{prop}
\begin{proof}
By proposition \ref{PropSdPFactorisesdP}, the adjunction
\begin{equation*}
\RealNP{\sd_P(-)}\colon \sS_P\leftrightarrow\Top_{N(P)}\colon \Ex_P\Sing_{N(P)}
\end{equation*}
is isomorphic to the adjunction $(\RealNP{C(\Sd_P(-))},\ex_P D \Sing_{N(P)})$, which can be rewritten as $(\C\Sd_P,\ex_P\D)$. By \cite[Theorem 2.12]{TopNP} and Theorem \ref{TheoCDQuillenEquivalence}, this is a composition of Quillen equivalences, so it must be a Quillen equivalence. But then, since all objects of $\Top_{N(P)}$ are fibrant and all objects of $\sS_P$ are cofibrant, the co-unit 
\begin{equation*}
\RealNP{\sd_P(\Ex_P\Sing_{N(P)}\fil{X})}\to \fil{X}
\end{equation*}
is a weak-equivalence in $\Top_{N(P)}$. By Proposition \ref{PropRealSdCW} the domain of this map admits the structure of a vertical filtered CW-complex. Furthermore, since both objects are fibrant and cofibrant, the map must be a filtered homotopy equivalence.
\end{proof}

\begin{corollaire}\label{CorollaireCategorieHomotopieCW}
Let $\Ver_P/{\sim}$ be the category whose objects are vertical filtered CW-complexes and whose  maps are classes of strongly filtered maps up to filtered homotopy. Then there is an equivalence of categories 
\begin{equation*}
\Ver_P/{\sim}\  \simeq \Ho(Top_{N(P)})
\end{equation*}
\end{corollaire}

\begin{proof}
Since $\Top_{N(P)}$ is a model category and all of its objects are fibrant, there is an equivalence of categories $\text{Cof}/{\sim}\ \simeq \Ho(\Top_{N(P)})$, where $\text{Cof}\subset \Top_{N(P)}$ is the full subcategory of cofibrant objects. But, since by proposition \ref{PropCofibrantVerticalCW}, all cofibrant objects of $\Top_{N(P)}$ are filtered homotopy equivalent to vertical filtered CW-complexes, one can replace $\text{Cof}$ by $\Ver_P$.
\end{proof}
\subsection{Proving the Quillen equivalence between $\Top_P$ and $\Top_{N(P)}$}
\label{SubsectionProvingQETopNP}

In this subsection, we prove the following.
\begin{theo}\label{TheoremeAMontrer}
Let $\fil{X}$ be a cofibrant object in $\Top_{N(P)}$, the unit
\begin{equation}\label{EquationAMontrer}
\fil{X}\to(X,\varphi_P\circ\varphi_X)\times_P\Real{N(P)}
\end{equation}
is a weak-equivalence in $\Top_{N(P)}$.
\end{theo}
Since the adjunction between $\Top_P$ and $\Top_{N(P)}$ preserves filtered homotopy equivalences, by proposition \ref{PropCofibrantVerticalCW}, it is enough to show the result for vertical filtered CW-complexes. We first make two observations :

\begin{lemme}\label{LemmeNPP}
The map $\RealNP{N(P)}\to \RealP{N(P)}\times_P\Real{N(P)}$ is a weak-equivalence in $\Top_{N(P)}$.
\end{lemme}

\begin{proof}
by \cite[Lemma 2.16]{TopNP} the map $\RealP{N(P)}\to P$ is a trivial fibration in $\Top_P$. Since the functor $-\times_P\times\Real{N(P)}$ preserves trivial fibrations, the map
\begin{equation*}
\RealP{N(P)}\times_P\Real{N(P)}\to P\times_P\Real{N(P)}\simeq \Real{N(P)}
\end{equation*}
is a trivial fibration in $\Top_{N(P)}$. Furhermore, the composition
\begin{equation*}
\Real{N(P)}\to \RealP{N(P)}\times_P\Real{N(P)}\to  \Real{N(P)}
\end{equation*}
is the identity. By two out of three, the first map is a weak-equivalence in $\Top_{N(P)}$.
\end{proof}

\begin{lemme}\label{LemmePullback}
Let $\fil{X}$ be an object in $\Top_{N(P)}$, the following commutative diagram is a pull-back square :
\begin{equation*}
\begin{tikzcd}
\fil{X}
\arrow{r}{(\Id_X,\varphi_X)}
\arrow{d}{\varphi_X}
&\filP{X}\times_P\Real{N(P)}
\arrow{d}{(\varphi_X,\Id)}
\\
\RealNP{N(P)}
\arrow{r}{(\Id,\Id)}
&\RealP{N(P)}\times_P\Real{N(P)}
\end{tikzcd}
\end{equation*}
\end{lemme}
\begin{proof}
This follows from a direct calculation of the pull-back.
\end{proof}
Those two observations imply the following lemma. 
\begin{lemme}\label{LemmeProductCase}
Let $X$ be some topological space and $\Delta^{\varphi}$ a non-degenerate filtered simplex. Then, the map
\begin{equation*}
X\otimes\RealNP{\Delta^{\varphi}}\to \left(X\otimes\RealP{\Delta^{\varphi}}\right)\times_P\Real{N(P)}
\end{equation*}
is a weak-equivalence in $\Top_{N(P)}$.
\end{lemme}

\begin{proof}
Write $\Delta^{\varphi}=[p_0,\dots,p_n]$, and write $Q=\{p_0<\dots<p_n\}\subset P$ for the corresponding linear sub-poset. Then $\Real{\Delta^{\varphi}}=\Real{N(Q)}$. Furthermore, by lemma \ref{LemmePullback} we have the following pullback square in the category $\Top_{N(Q)}$
\begin{equation*}
\begin{tikzcd}
X\otimes||\Delta^{\varphi}||_{N(Q)}
\arrow{r}
\arrow{d}{\pr_{\Delta^{\varphi}}}
&\left(X\otimes||\Delta^{\varphi}||_{Q}\right)\times_{Q}\Real{N(Q)}
\arrow{d}
\\
\Real{N(Q)}_{N(Q)}
\arrow{r}{(\Id,\Id)}
&\Real{N(Q)}_{Q}\times_{Q}\Real{N(Q)}
\end{tikzcd}
\end{equation*}
We know by lemma \ref{LemmeNPP} that the bottom map is a weak-equivalence in $\Top_{N(Q)}$. Since $\Top_{N(Q)}$ is a right proper model category, it is enough to show that the right map is a fibration in $\Top_{N(Q)}$ to deduce that the top map is a weak-equivalence in $\Top_{N(Q)}$. But the map $\pr_{\Delta^{\varphi}}\colon X\otimes\Real{\Delta^{\varphi}}_Q\to \Real{N(Q)}_Q$ is clearly a fibration in $\Top_Q$, since it is a projection. In turn, its image by the right Quillen functor $-\times_Q\Real{N(Q)}$ is a fibration in $\Top_{N(Q)}$. Now, the inclusion $Q\subset P$ induces a functor $\Top_{N(Q)}\to \Top_{N(P)}$ that preserves all weak-equivalences, and that sends the map $X\otimes||\Delta^{\varphi}||_{N(Q)}\to \left(X\otimes||\Delta^{\varphi}||_{Q}\right)\times_{Q}\Real{N(Q)}$ to the map $X\otimes||\Delta^{\varphi}||_{N(P)}\to \left(X\otimes||\Delta^{\varphi}||_{P}\right)\times_{P}\Real{\Delta^{\varphi}}$. So it remains to be shown that the inclusion
\begin{equation*}
\left(X\otimes||\Delta^{\varphi}||_{P}\right)\times_{P}\Real{\Delta^{\varphi}}\to \left(X\otimes||\Delta^{\varphi}||_{P}\right)\times_{P}\Real{N(P)}
\end{equation*}
is a weak-equivalence in $\Top_{N(P)}$. Let $\Delta^{\psi}$ be a non-degenerate filtered simplex and consider an element, $f$, in $\Map(\RealNP{\Delta^{\psi}},\left(X\otimes||\Delta^{\varphi}||_{P}\right)\times_{P}\Real{N(P)})$, of dimension $n$. It can be decomposed as the product of a map in $\Top_P$, $f_P\colon \Delta^n\otimes\Delta^{\psi}\to X\otimes\RealP{\Delta^{\varphi}}$, and the inclusion $\RealNP{\Delta^{\psi}}\to \RealNP{N(P)}$. In particular, for $f_P$ to exist, we must have $\Delta^{\psi}\subset \Delta^{\varphi}$. We conclude that $f$ must lie in $\Map(\RealNP{\Delta^{\psi}},\left(X\otimes||\Delta^{\varphi}||_{P}\right)\times_{P}\Real{\Delta^{\varphi}})$, which conclude the proof.

\end{proof}

Before moving on to the proof of Theorem \ref{TheoremeAMontrer}, we need a few preliminary results and definitions

\begin{lemme}\label{LemmeWEIffLift}
Let $\alpha\colon \fil{X}\to\fil{Y}$ be a map in $\Top_{N(P)}$. It is a weak-equivalence if and only if, for all commutative square of the form
\begin{equation}\label{EquationLiftingWeakEquivalence}
\begin{tikzcd}
\fil{X}
\arrow{r}{\alpha}
&\fil{Y}
\\
S^{n-1}\otimes\Delta^{\varphi}
\arrow{r}
\arrow{u}{f}
&B^n\otimes\Delta^{\varphi}
\arrow{u}{g}
\arrow[dashed]{ul}{h}
\end{tikzcd}
\end{equation}
where $n\geq 0$ and $\Delta^{\varphi}$ is a non-degenerate filtered simplex, there exists $h\colon B^n\otimes\Delta^{\varphi}\to\fil{X}$ such that $h_{|S^{n-1}\otimes\Delta^{\varphi}}=f$ and $\alpha\circ h$ is filtered homotopic to $g$ relative to $S^{n-1}\otimes\Delta^{\varphi}$.

The result also holds in $\Top_P$.
\end{lemme}

\begin{proof}
To prove the direct implication, consider a factorisation of $\alpha$ as follows
\begin{equation*}
\begin{tikzcd}
\fil{X} \arrow[bend right]{r}{j}
& \fil{Z}
\arrow[bend right]{l}{r}
\arrow{r}{p}
&\fil{Y}
\end{tikzcd}
\end{equation*}
such that, $\alpha = p\circ j$, with $p$ a fibration, $r\circ j= \Id_X$ and $j\circ r$ is filtered homotopic to $\Id_Z$, relative to $j(X)$. Such a decomposition can be produced by considering the path-space associated to $\alpha$ (see \cite[Lemma 2.11]{TopNP}).
We now have the following lifting problem
\begin{equation*}
\begin{tikzcd}
\fil{X} \arrow{r}{j}
& \fil{Z}
\arrow{r}{p}
&\fil{Y}
\\
S^{n-1}\otimes\Delta^{\varphi}
\arrow{rr}
\arrow{u}{f}
&&B^n\otimes\Delta^{\varphi}
\arrow{u}{g}
\arrow[dashed]{ul}{h'}
\end{tikzcd}
\end{equation*}
Since $\alpha$ is a weak-equivalence by hypothesis, and since $j$ is a filtered homotopy equivalence, $p$ must be a trivial fibration. This means that there exists some lift $h'\colon \Delta^n\otimes\Delta^{\varphi}\to \fil{Z}$. 
Taking $h=r\circ h'$ gives the desired lift. 
Indeed, we have $r\circ h'_{|S^{n-1}\otimes\Delta^{\varphi}}=r\circ j\circ f=f$. On the other hand, $\alpha\circ h = \alpha\circ r\circ h' = p\circ j\circ r\circ h'$, which is filtered homotopic to $p\circ h' = g$. 
The homotopy between $j\circ r$ and $\Id_Z$ is constant on $j(X)$, in which lies $j\circ f( S^{n-1}\otimes\Delta^{\varphi})$, this implies that the homotopy between $\alpha\circ h$ and $g$ is relative to $S^{n-1}\otimes\Delta^{\varphi}$.

To prove the converse, we need to prove that $f$ induces isomorphisms on all filtered homotopy groups. Let $\phi\colon \Delta^{\varphi}\to \fil{X}$ be some pointing of $\fil{X}$. Any element in $\pi_n(\Map(\Delta^{\varphi},\fil{X}),\phi)$ can be represented by a map $f\colon S^n\otimes\Delta^{\varphi}\to \fil{X}$, sending $\{*\}\otimes\Delta^{\varphi}$ to the chosen pointing. Assume that the element represented by $\alpha\circ f $ is trivial in $\pi_n(\Map(\Delta^{\varphi},\fil{Y}),\alpha\circ\phi)$. This means that $\alpha\circ f$ extends to a map $g\colon B^{n+1}\otimes\Delta^{\varphi}\to \fil{Y}$. Now, by the lifting property \eqref{EquationLiftingWeakEquivalence}, we deduce that there exists a map $h\colon B^{n+1}\otimes\Delta^{\varphi}\to \fil{X}$ extending $f$, which means that $f$ represents the trivial element in $\pi_n(\Map(\Delta^{\varphi},\fil{X}),\phi)$. This implies that $\alpha$ induces an injective map on filtered homotopy groups. To show the surjectivity, assume that $g\colon S^n\otimes \Delta^{\varphi}\to \fil{Y}$ is a representant of some element in $\pi_n(\Map(\Delta^{\varphi},\fil{Y}),\alpha\circ \phi)$. The map $g$ can be seen as a map $\widetilde{g}\colon B^n\otimes\Delta^{\varphi}\to \fil{Y}$ sending $S^{n-1}\otimes\Delta^{\varphi}$ to $\Im(\alpha\circ\phi)\simeq \Delta^{\varphi}$. 
In particular, the restriction of $\widetilde{g}$ to $S^{n-1}\otimes\Delta^{\varphi}$ lifts to $\fil{X}$, where it lands in $\Im(\phi)$. The lifting property \eqref{EquationLiftingWeakEquivalence} then guarantees that there exists some map $h\colon B^n\otimes\Delta^{\varphi}\to \fil{X}$, sending $S^{n-1}\otimes\Delta^{\varphi}$ to the pointing of $X$, whose image by $\alpha$ is filtered homotopic to $\widetilde{g}$ relative to $S^{n-1}\otimes\Delta^{\varphi}$. In particular, the lift $h$ represents an element of $\pi_n(\Map(\Delta^{\varphi},\fil{X}),\phi)$ in the preimage of the element represented by $g$. This implies that $\alpha$ induces a surjection on filtered homotopy groups.
The same proof works in the case of $\Top_P$.
\end{proof}

\begin{construction}
Let $\fil{X}$ be a vertical filtered CW-complex, and $\Delta^{\varphi}$ a non-degenerate filtered simplex. Let $L(X)^{\varphi}$ be the subcomplex of $L(X)$ containing all cells with label $\Delta^{\psi}$ such that $\Delta^{\varphi}\subset\Delta^{\psi}$. Let $\lambda^{\varphi}$ be the constant labeling map with value $\Delta^{\varphi}$. We define $\fil{X}^{\varphi}=V(L(X)^{\varphi},\lambda^{\varphi})$. 
Note that the we have a well defined composition of label-preserving maps
\begin{equation*}
\begin{tikzcd}
(L(X)^{\varphi},\lambda^{\varphi})
\arrow{r}{\Id_{L(X)^{\varphi}}} 
&(L(X)^{\varphi},\lambda_{L(X)})
\arrow[hookrightarrow]{r}
&(L(X),\lambda_{L(X)})
\end{tikzcd}
\end{equation*}
Applying $V$ to this map gives a well-defined monomorphism $\fil{X}^{\varphi}\hookrightarrow \fil{X}$, though it is not the inclusion of a subcomplex. By abuse of notation, we will write $X^{\varphi}$ for the subspace of $X$ underlying $\fil{X}^{\varphi}$.
\end{construction}

\begin{lemme}\label{LemmeFLandsInXPhi}
Let $f\colon K\otimes\Delta^{\varphi}\to \fil{X}$ be a map in $\Top_{N(P)}$ with $K$ some topological space, $\Delta^{\varphi}$ a non-degenerate filtered simplex and $\fil{X}$ a vertical filtered CW-complex. Then, $\Im(f)\subset X^{\varphi}$.
\end{lemme}

\begin{proof}
By assumption, $f$ is strongly filtered. This implies that for any $x\in K$, and any $t$ in the interior of $\Delta^{\varphi}$, $f(x,t)$ must lie in a cell of dimension $(n,\Delta^{\psi})$, with $\Delta^{\varphi}\subset \Delta^{\psi}$. By continuity of $f$, this implies that $\Im(f)$ must lie in the closure of the union of all cells of dimension $(n,\Delta^{\psi})$ with $\Delta^{\varphi}\subset \Delta^{\psi}$. Any cell (say of dimension $(m,\Delta^{\mu})$) intersecting this closure must intersect the boundary of some cell of dimension $(n,\Delta^{\psi})$ with $\Delta^{\varphi}\subset\Delta^{\psi}$. But this implies that $\Delta^{\psi}\subset\Delta^{\mu}$. Now, if $f(x,t)$ lies in some cell of dimension $(n,\Delta^{\psi})$, the fact that $f$ is strongly filtered implies that $f(x,t)=(y,t)\in \Delta^n\otimes\Delta^{\varphi}\subset \Delta^n\otimes\Delta^{\psi}$ for some $y\in \Delta^n$. In particular, $f(x,t)$ must lie in $\fil{X}^{\varphi}$.
\end{proof}

\begin{remarque}
Lemma \ref{LemmeFLandsInXPhi} implies that the stratified homotopy groups of a vertical filtered CW-complex can be computed from its associated P-labeled CW-complex. Indeed, let $\fil{X}$ be some vertical filtered CW-complex, $\Delta^{\varphi}$ be a non-degenerate filtered simplex, and $\phi\colon \Real{V}\to \fil{X}$ a pointing of $\fil{X}$. For simplicity, we will assume that $\phi$ is the inclusion of some $0$-cell. Then, any element in $s\pi_n(\fil{X},\phi)(\Delta^{\varphi})$ can be represented by a map $f\colon S^n\otimes\Delta^{\varphi}\to \fil{X}$. But, by lemma \ref{LemmeFLandsInXPhi}, such a map must land in $X^{\varphi}\simeq L(X)^{\varphi}\otimes\Delta^{\varphi}$. In particular, one can define $\widetilde{f}\colon S^n\otimes\Delta^{\varphi}\to L(X)^{\varphi}\otimes\Delta^{\varphi}$ by $\widetilde{f}(x,(t_0,\dots,t_n))=(\pr_1(f(x,(1,0,\dots,0))),(t_0,\dots,t_n))$. The map $\widetilde{f}$ is filtered homotopic to $f$, by construction, and is of the form $V(g\colon S^n\to L(X)^{\varphi})$ for some map $g$ between $P$-labeled CW-complexes. Doing the same for maps of the form $B^{n+1}\otimes\Delta^{\varphi}\to \fil{X}$, we get an isomorphism
\begin{equation*}
s\pi_n(\fil{X},\phi)(\Delta^{\varphi})\simeq \pi_n(L(X)^{\varphi},\phi).
\end{equation*}

\end{remarque}

\begin{construction}
Let $\fil{X}$ be a vertical filtered CW-complex, and let $p\in P$. 
Let $L(X)^{\leq p}$ be the subcomplex of $L(X)$ containing all cells with label $\Delta^{\psi}$ such that $p\in \Delta^{\psi}$. For $\Delta^{\psi}=[q_0,\dots,q_m]$ a non-degenerate filtered simplex, such that $p=q_i$ for some $0\leq i \leq m$, define $\tr_{\leq p}(\Delta^{\psi})=[q_0,\dots,q_i]$. We define $\fil{X}^{\leq p}=V(L(X)^{\leq p },\tr_{\leq p }\circ \lambda_{L(X)})$. As before, we will write $X^{\leq p}$ to denote the subspace of $X$ underlying $\fil{X}^{\leq p }$.

Now, let $\Delta^{\varphi}=[p_0,\dots,p_n]$ be a non-degenerate filtered simplex. Define the subspace $X^{\leq \varphi}\subset X^{\leq p_n}$ to be the following union :
\begin{equation*}
X^{\leq \varphi}= \bigcup\limits_{p_i\in \Delta^{\varphi}} (\varphi_P\circ\varphi_X)^{-1}({p_i})\cap X^{\leq p_n}
\end{equation*}
In other word $X^{\leq \varphi}$ is the union of the $p_i$-strata of $X^{\leq p_n}$.
\end{construction}

\begin{lemme}\label{LemmeGLandsInXleqPhi}
Let $g\colon K\otimes\Delta^{\varphi}\to \filP{X}$ be a map in $\Top_P$, with $K$ some topological space, $\Delta^{\varphi}$ a non-degenerate filtered simplex and $\fil{X}$ a vertical filtered CW-complex. Then $\Im(g)$ lies in $X^{\leq\varphi}$.
\end{lemme}

\begin{proof}
Consider the following subset of $\Real{\Delta^{\varphi}}=\Real{[p_0,\dots,p_n]}$ 
\begin{equation*}
A=\{(t_0,\dots,t_n)\ | \ \sum t_i =1, \ 0\leq t_i \leq 1,\ \forall\ 0\leq i\leq n-1, \ 0<t_n\leq 1\}.
\end{equation*}
it is dense in $\Real{\Delta^{\varphi}}$, and all points in $A$ are mapped to $p_n$ by the stratification $\varphi_P$. Since $g$ is a map in $\Top_P$, it preserves the stratifications over $P$, and so must send $K\otimes A$ to the $p_n$-stratum of $\filP{X}$. But the $p_n$-stratum of $\filP{X}$ is contained in $X^{\leq p_n}$, which is a closed subset of $X$. This implies that $\Im(g)$ lies in $X^{\leq p_n}$. Furthermore, any point in $\Im(g)$ lies in the $p_i$-stratum of $\filP{X}$ for some $0\leq i\leq n$, which implies that $\Im(g)$ lies in $X^{\leq\varphi}$.
\end{proof}

\begin{lemme}\label{LemmeUvarphiDeformationRetracts}
Let $\fil{X}$ be a vertical filtered CW-complex, and $\Delta^{\varphi}$ be a non-degenerate filtered simplex. There exists some open $U_{\varphi}\subset X^{\leq \varphi}$ such that $X^{\varphi}\subset U_{\varphi}$, and $(U_{\varphi},\varphi_P\circ\varphi_X)$ deformation retracts to $(X^{\varphi},\varphi_P\circ\varphi_X)$ in the category $\Top_P$.
\end{lemme}

\begin{proof}
Consider the inclusion of CW-complexes $L(X)^{\varphi}\subset L(X)^{\leq p}$. By classical results, \cite[Proposition A.5]{Hatcher} There exists an open neighborhood of $L(X)^{\varphi}$ into $L(X)^{\leq p}$, $W_{\varphi}$, that deformation retracts to $L(X)^{\varphi}$. Let us write $i\colon L(X)^{\varphi}\to W_{\varphi}$ for the inclusion, $r\colon W_{\varphi}\to L(X)^{\varphi}$ for the retraction and $H\colon W_{\varphi}\times [0,1]\to W_{\varphi}$ for the homotopy. 
Note that $H$ can be chosen such that for a point $x$ in the interior of a cell $e_{\alpha}$, $H(x,t)$ lies in the closure of $e_{\alpha}$ for all $t$. In particular, those maps are label-preserving. 
Now define a $P$-labeled CW-complex $\lab{Z}$ as follows. Set $Z= L(X)^{\leq p}$, and $\lambda_Z(e_{\alpha})=\lambda_{L(X)}(e_{\alpha})\cap \Delta^{\varphi}$. This is well defined since for any cell $e_{\alpha}\in L(X)^{\leq p}$, $p\in \lambda_{L(X)}(e_{\alpha})\cap \Delta^{\varphi}$. Now, define $\fil{Y}=V\lab{Z}$. Clearly, there is a sequence of inclusions of subspaces
\begin{equation*}
\fil{X}^{\varphi}\subset \fil{Y}\subset \fil{X}^{\leq \varphi}
\end{equation*}
Note that since $\fil{Y}$ is a filtered subspace of $\fil{X}^{\leq p}$, $\varphi_Y$ is simply the restriction of $\varphi_X$ to $Y$.
Furthermore, the open subset, $(W_{\varphi},\lambda_{L(X)})\subset (L(X)^{\leq p},\lambda_{L(X)}),$ lifts to some open subset $V(W_{\varphi})\subset X^{\leq p}$. Define $U_{\varphi}= V(W_{\varphi})\cap X^{\leq \varphi}$. We then have the sequence of inclusions
\begin{equation*}
\begin{tikzcd}
(X^{\varphi},\varphi_P\circ\varphi_X)
\arrow[hookrightarrow]{r}{j_1}
&(Y\cap U_{\varphi},\varphi_P\circ \varphi_X)
\arrow[hookrightarrow]{r}{j_2}
&(U_\varphi,\varphi_P\circ \varphi_X)
\end{tikzcd}
\end{equation*}
We will show that both admit deformation retracts. First, consider some cell $e_{\alpha}$ of $Y$ such that $e_{\alpha}$ intersects $U_{\varphi}$, and $e_{\alpha}$ is not a cell of $X^{\varphi}$. By construction of $Y$, there must be some cell $e'_{\alpha}\in Z$ such that $e_{\alpha}=V(e'_{\alpha})\simeq e'_{\alpha}\otimes \lambda_Z(e'_{\alpha})$. 
But then, the intersection $e_{\alpha}\cap U_{\varphi}$ is filtered homeomorphic to $\left(e'_{\alpha}\cap W_{\varphi}\right)\otimes \lambda_Z(e'_{\alpha})$. 
In particular, using this homeomorphism, we can define $r_{1,\alpha}\colon e_{\alpha}\cap U_{\varphi}\to X^{\varphi}$ as a product of $r_{\alpha}\colon e'_{\alpha}\cap W_{\varphi}\to L(X)^{\varphi}$ and the inclusion $\lambda_Z(e'_{\alpha})\subset \Delta^{\varphi}$, 
and the homotopy $H_{1,\alpha}\colon (e_{\alpha}\cap U_{\varphi})\times [0,1]\to \left(\overline{e_{\alpha}}\cap U_{\varphi}\right)$ as the product of $H$ and the identity of $\lambda_Z(e'_{\alpha})$. 
The continuity of the maps $H$ and $r$ then guarantees that the $r_{1,\alpha}$ and $H_{1,\alpha}$ glue together to form a retraction $r_1\colon Y\cap U_{\varphi}\to X^{\varphi}$ and a homotopy $H_{1}\colon Y\cap U_{\varphi}\times [0,1]\to Y\cap U_{\varphi}$ between the identity of $Y\cap U_{\varphi}$ and $r_1\circ j_1$. One note that all the maps involved are in fact vertical. In particular, $j_1$ admits a deformation retraction in the category $\Top_{N(P)}$, which means it also does in the category $\Top_P$.

Now, let $e_{\alpha}=V(e'_{\alpha})$ be a cell in $X^{\leq p}$ of dimension $(k,\Delta^{\psi_{\alpha}})$. The intersection $e_{\alpha}\cap U_{\varphi}$ can be further decomposed as $(e_{\alpha}\cap V(W_{\varphi}))\cap X^{\leq \varphi}$. In particular, there is a filtered homeomorphism in $\Top_P$
\begin{equation*}
e_{\alpha}\cap U_{\varphi} \simeq \left(e'_{\alpha}\cap W_{\varphi}\right)\otimes \RealP{\Delta^{\psi_{\alpha}}}^{\leq \varphi}
\end{equation*}
Where, $\RealP{\Delta^{\psi_{\alpha}}}^{\leq \varphi}$ is the union of the $p_i$-strata of $\RealP{\Delta^{\psi_{\alpha}}}$. Alternatively, if we write $\Delta^{\psi_{\alpha}}=[q_0,\dots,q_m]$, we can describe the subset $\RealP{\Delta^{\psi_{\alpha}}}^{\leq \varphi}\subset \RealP{\Delta^{\psi_{\alpha}}}$ as follows :
\begin{equation*}
\RealP{\Delta^{\psi_{\alpha}}}^{\leq \varphi}=\{(t_0,\dots,t_m)\ | \max\{q_i\ | \ t_i\not= 0\}\in \Delta^{\varphi}\}
\end{equation*}
One can then define a filtered retract $\RealP{\Delta^{\psi_{\alpha}}}^{\leq \varphi}\to \RealP{\Delta^{\psi_{\alpha}}\cap \Delta^{\varphi}}$ as follows. 
If $(t_0,\dots,t_m)\in \RealP{\Delta^{\psi_{\alpha}}}^{\leq \varphi}$, define $(t_0,\dots,t_m)_{\varphi}=(s_0,\dots,s_m)$ with $s_i=t_i$ if $q_i\in \Delta^{\varphi}$, and $s_i=0$ else, and write $|(s_0,\dots,s_m)|=\sum_i s_i$. Then, one can define the retraction as follows
\begin{align*}
r_{2,\psi_{\alpha}}\colon \RealP{\Delta^{\psi_{\alpha}}}^{\leq \varphi}&\to \RealP{\Delta^{\psi_{\alpha}}\cap\Delta^{\varphi}}\\
(t_0,\dots,t_m)&\mapsto \frac{(t_0,\dots,t_m)_{\varphi}}{|(t_0,\dots,t_m)_{\varphi}|}
\end{align*}
We can then take the linear homotopy, $H_{2,\alpha}$, between the identity of $\RealP{\Delta^{\psi_{\alpha}}}^{\leq\varphi}$ and the composition
\begin{equation*}
\begin{tikzcd}
\RealP{\Delta^{\psi_{\alpha}}}^{\leq\varphi}
\arrow{r}{r_{2,\alpha}}
&\RealP{\Delta^{\psi_{\alpha}}\cap\Delta^{\varphi}}
\arrow{r}{j_{2,\alpha}}
&\RealP{\Delta^{\psi_{\alpha}}}^{\leq\varphi}
\end{tikzcd}
\end{equation*}
Note that by construction, the maps $j_{2,\alpha}$, $r_{2,\alpha}$ and $H_{2,\alpha}$ are filtered over $P$, and are compatible with inclusions $\Delta^{\psi_{\alpha}}\hookrightarrow \Delta^{\psi_{\beta}}$. Furthermore, taking the product of $j_{2,\alpha}$ with the identity of $e'_{\alpha}$, we get a map (that we will still write $j_{2,\alpha}$) 
\begin{equation*}
j_{2,\alpha}\colon \left(e'_{\alpha}\cap W_{\varphi}\right)\otimes \RealP{\Delta^{\psi_{\alpha}}\cap \Delta^{\varphi}}\to \left(e'_{\alpha}\cap W_{\varphi}\right)\otimes\RealP{\Delta^{\psi_{\alpha}}}^{\leq\varphi}
\end{equation*}
In particular, the $j_{2,\alpha}$ can be glued together to produce the inclusion $j_2\colon Y\cap U_{\varphi}\to U_{\varphi}$. Doing the same for the $r_{2,\alpha}$ produces a retract $r_2\colon U_{\varphi}\to Y\cap U_{\varphi}$ and similarly for the $H_{2,\alpha}$ giving a homotopy between $j_2\circ r_2$ and $\Id_{U_{\varphi}}$. Since the $r_{2,\alpha}$ and $H_{2,\alpha}$ are filtered maps over $P$, so are $r_2$ and $H_2$. In particular, $j_2$ admits a deformation retraction in $\Top_P$.
\end{proof}

\begin{proof}[Proof of Theorem \ref{TheoremeAMontrer}]
Let $\fil{X}$ be a filtered CW-complexes. By Lemma \ref{LemmeWEIffLift}, it is enough to show that for any $n\geq 0$, any non-degenerate filtered simplex $\Delta^{\varphi}$, and any commutative diagram in $\Top_{N(P)}$ of the form
\begin{equation}\label{EquationLiftingProblemBeforeDecomposition}
\begin{tikzcd}
\fil{X}
\arrow{r}{\alpha}
&(X,\varphi_P\circ\varphi_X)\times_ P \Real{N(P)}
\\
S^{n-1}\otimes\Delta^{\varphi}
\arrow{r}
\arrow{u}{f}
&B^n\otimes\Delta^{\varphi}
\arrow[swap]{u}{g}
\arrow[dashed]{ul}{h}
\end{tikzcd}
\end{equation}
there exists a lift up to homotopy $h\colon B^n\otimes\Delta^{\varphi}\to \fil{X}$ in $\Top_{N(P)}$. By lemma \ref{LemmeFLandsInXPhi}, we can replace $\fil{X}$ with the subspace $\fil{X}^{\varphi}$. By lemma \ref{LemmeGLandsInXleqPhi}, we can replace $\filP{X}\times_P\Real{N(P)}$ by $(X^{\leq\varphi},\varphi_P\circ\varphi_X)\times_P\Real{N(P)}$. The lifting problem \eqref{EquationLiftingProblemBeforeDecomposition} then turns into 
\begin{equation*}
\begin{tikzcd}
\fil{X}^{\varphi}
\arrow{r} 
&(X^{\leq\varphi},\varphi_P\circ\varphi_X)\times_ P \Real{N(P)}
\\
S^{n-1}\otimes\Delta^{\varphi}
\arrow{r}
\arrow{u}{f}
&B^n\otimes\Delta^{\varphi}
\arrow[swap]{u}{g}
\arrow[dashed]{ul}{h}
\end{tikzcd}
\end{equation*}
We then decompose the lifting problem as follows (note that in future diagrams, we omit the stratifications).
\begin{equation*}
\begin{tikzcd}[row sep = 80]
X^{\varphi}
\arrow{r}{\alpha}
&X^{\varphi}\times_P\Real{N(P)}
\arrow[hookrightarrow]{r}
&U_{\varphi}\times_P\Real{N(P)}
\arrow[hookrightarrow]{r}
&X^{\leq\varphi}\times_ P \Real{N(P)}
\\
S^{n-1}\otimes\Delta^{\varphi}
\arrow{rrr}
\arrow{u}{f}
&&&B^n\otimes\Delta^{\varphi}
\arrow[swap]{u}{g}
\arrow[dashed]{ulll}{h_3}
\arrow[dashed,swap]{ull}{h_2}
\arrow[dashed,swap]{ul}{h_1}
\end{tikzcd}
\end{equation*}
We will construct the $h_i$ sequentially. Note that while working on the right side of the lifting problem, we can omit the "$\times_P\Real{N(P)}$" and consider that all spaces and maps live in $\Top_P$. Consider the following lifting problem.
\begin{equation*}
\begin{tikzcd}
U_{\varphi}
\arrow[hookrightarrow]{r}
&X^{\leq\varphi}
\\
S^{n-1}\otimes\Delta^{\varphi}
\arrow{r}
\arrow{u}{\alpha\circ f}
&B^n\otimes\Delta^{\varphi}
\arrow[swap]{u}{g}
\arrow[dashed,swap]{ul}{h_1}
\end{tikzcd}
\end{equation*}
Write $\Delta^{\varphi}=[p_0,\dots,p_n]$. Let us show that $g(B^n\otimes\{p_0\})$ lies in $U_{\varphi}$. If $\Delta^{\varphi}=\{p_0\}$, this follows from the definition of $U_{\varphi}$, and tautologically implies the existence of a lift $h_1$. 
For the other cases, let $(x,(1,0,\dots,0))$ be a point in $B^n\otimes\{p_0\}$, the point $g(x,(1,0,\dots,0))$ must lie in some cell $e_{\alpha}$ of $X^{\leq p}$, of dimension $(k,\Delta^{\psi_{\alpha}})$. 
But for any $p_i\in \Delta^{\varphi}$, $g(x,(1,0,\dots,0))$ must lie in the closure of the $p_i$-stratum of $X^{\leq p}$.
Indeed, consider the sequence $g(x,(1-\frac{1}{l},0,\dots,0,\frac{1}{l},0,\dots,0))$, where $\frac{1}{l}$ is in the $(i+1)$-th position. It lies in the $p_i$-stratum of $X$ but converges to $g(x,(1,0,\dots,0))$ when $l$ tends to infinity. 
This implies that $p_i\in \Delta^{\psi_{\alpha}}$ for all $p_i\in \Delta^{\varphi}$, which means that $\Delta^{\varphi}\subset\Delta^{\psi_{\alpha}}$. 
This implies that $W_{\varphi}$ contains the (non-filtered) cell $e'_{\alpha}$ which satisfies $V(e'_{\alpha})=e_{\alpha}$. In turn, $V(W_{\varphi})$ contains the cell $e_{\alpha}$, and so the $p_0$-stratum of $e_{\alpha}$ is contained in $U_{\varphi}$. In particular, $U_{\varphi}$ must contain the image of $B^n\otimes\{p_0\}$ under $g$. Now, consider the following subsets of $B^n\otimes\Delta^{\varphi}$, for $\epsilon>0$.
\begin{equation*}
\tr_{\epsilon}(B^n\otimes\Delta^{\varphi})=\{(x,(t_0,\dots,t_n))\ | t_0 \geq  1-\epsilon\}\subset B^n\otimes\Delta^{\varphi}
\end{equation*}
For all $\epsilon>0$, $\tr_{\epsilon}(B^n\otimes\Delta^{\varphi})$ is compact, which means that $(g(\tr_{\epsilon}(B^n\otimes\Delta^{\varphi})))_{\epsilon >0}$ is a family of nested compact in $X^{\leq \varphi}$. 
As we have shown, their intersection lies in the open set $U_{\varphi}$, which implies that their exists some $\epsilon>0$ such that $g(\tr_{\epsilon}(B^n\otimes\Delta^{\varphi}))$ lies in $U_{\varphi}$. 
On the other hand, fix some homeomorphism between $S^{n-1}\times[0,1[$ and $B^n\setminus \{0\}$. The family $(g(S^{n-1}\times[0,\delta]\otimes\Delta^{\varphi}))_{0<\delta <1}$, is a family of nested compacts in $X^{\leq \varphi}$. Their intersection is $g(S^{n-1}\times\{0\}\otimes\Delta^{\varphi})=\alpha\circ f(S^{n-1}\otimes\Delta^{\varphi})\subset U_{\varphi}$. This implies that their exist $0<\delta < 1$ such that $g(S^{n-1}\times[0,\delta]\otimes\Delta^{\varphi})$ lies in $U_{\varphi}$. 
Let $A^{n,\varphi}(\epsilon,\delta)=\tr_{\epsilon}(B^n\otimes\Delta^{\varphi})\cup S^{n-1}\times[0,\delta]\otimes\Delta^{\varphi} \subset B^n\otimes\Delta^{\varphi}$, and write $j$ for the corresponding inclusion. We have shown that $g(A^{n,\varphi}(\epsilon,\delta))$ lies in $U_{\varphi}$. 
By lemma \ref{LemmeAEpsilonPhi}, there exists a retraction $r\colon B^n\otimes\Delta^{\varphi}\to A^{n,\varphi}(\epsilon,\delta)$, and a filtered homotopy 
$H_1\colon (B^n\otimes\Delta^{\varphi})\times [0,1]\to B^n\otimes\Delta^{\varphi}$ between $\Id_{B^n\otimes\Delta^{\varphi}}$ and $j\circ r$, relative to $A^{n,\varphi}(\epsilon,\delta)$ in $\Top_P$. We then obtain the lift $h_1$ by taking $h_1=g\circ r$.

Now that $h_1$ has been constructed, consider the second commutative diagram in $\Top_P$ :
\begin{equation*}
\begin{tikzcd}
X^{\varphi}
\arrow[hookrightarrow]{r}
&U_{\varphi}
\\
S^{n-1}\otimes\Delta^{\varphi}
\arrow{u}{\alpha\circ f}
\arrow{r}
&B^n\otimes\Delta^{\varphi}
\arrow[swap]{u}{h_1}
\arrow[dashed,swap]{ul}{h_2}
\end{tikzcd}
\end{equation*}
By lemma \ref{LemmeUvarphiDeformationRetracts}, the top map is a filtered homotopy equivalence in $\Top_P$. In particular, it is a weak-equivalence in $\Top_P$. By lemma \ref{LemmeWEIffLift}, this implies that a lift up to homotopy, $h_2$, must exist.
Now turn to the last lifting problem :
\begin{equation*}
\begin{tikzcd}
(X,\varphi_X)^{\varphi}
\arrow{r}
&(X^{\varphi},\varphi_P\circ\varphi_X)\times_P\Real{N(P)}
\\
S^{n-1}\otimes\Delta^{\varphi}
\arrow{u}{f}
\arrow{r}
&B^n\otimes\Delta^{\varphi}
\arrow[swap]{u}{h_2}
\arrow[swap,dashed]{ul}{h_3}
\end{tikzcd}
\end{equation*}
By construction, $(X,\varphi_X)^{\varphi}\simeq L(X)^{\varphi}\otimes\Delta^{\varphi}$. By Lemma \ref{LemmeProductCase}, this means that the top map is a weak-equivalence in $\Top_{N(P)}$, which implies by lemma \ref{LemmeWEIffLift} that the desired lift must exist.
\end{proof}

\begin{lemme}\label{LemmeAEpsilonPhi}
The inclusion
\begin{equation*}
j\colon A^{n,\varphi}(\epsilon,\delta)=\tr_{\epsilon}(B^n\otimes\Delta^{\varphi})\cup (S^{n-1}\times [0,\delta])\otimes\Delta^{\varphi}\to B^n\otimes\Delta^{\varphi}
\end{equation*}
admit a deformation retraction $r\colon B^n\otimes\Delta^{\varphi}\to A^{n,\varphi}(\epsilon,\delta)$ in $\Top_P$.
\end{lemme}

\begin{proof}
Write $\Delta^{\varphi}=[p_0,\dots,p_k]$, and $\Delta^{\widehat{\varphi}}=[p_1,\dots,p_k]$. And fix an identification of  $B^n$ with $\{(x_1,\dots,x_n)\ |\ |(x_1,\dots,x_n)|\leq 1\}\subset \mathbb{R}^n$, such that the subspace $S^{n-1}\times [0,\delta]$ is identified with $\{(x_1,\dots,x_n)\ |\ 1-\delta\leq |(x_1,\dots,x_n)|\leq 1\}$.
Define the following subspaces of $B^n\otimes\Delta^{\varphi}$ :
\begin{align*}
B^{n,\varphi}(\epsilon,\delta)&=\{((x_1,\dots,x_n),(t_0,\dots,t_k))\ |\  |(x_1,\dots,x_n)|\leq 1-\delta, t_0\in [0,1-\epsilon]\}\\
C^{n,\varphi}(\epsilon,\delta)&=\{((x_1,\dots,x_n),(t_0,\dots,t_k))\ |\ \left( |(x_1,\dots,x_n)| = 1-\delta, \text{ and } t_0\in [0,1-\epsilon]\right)\\
&\phantom{\{((x_1,\dots,x_n),(t_0,\dots,t_k))\ |} \text{ or }\left( (x_0,\dots,x_n)| \leq 1-\delta, \text{ and } t_0= 1-\epsilon\right)\}
\end{align*}
By construction, we have
\begin{equation*}
B^n\otimes\Delta^{\varphi}=A^{n,\varphi}(\epsilon,\delta)\cup_{C^{n,\varphi}(\epsilon,\delta)}B^{n,\varphi}(\epsilon,\delta)
\end{equation*}
In particular, it is enough to show that the inclusion $C^{n,\varphi}(\epsilon,\delta)\hookrightarrow B^{n,\varphi}(\epsilon,\delta)$ admits a deformation retraction. Now, consider the following map :
\begin{align*}
\{(t_0,\dots,t_k)\ |t_0\in [0,1-\epsilon]\}&\mapsto [0,1]\times\RealP{\Delta^{\widehat{\varphi}}}\\
(t_0,\dots,t_k)&\mapsto (\frac{t_0}{1-\epsilon},\frac{(t_1,\dots,t_k)}{1-t_0})
\end{align*}
it is a filtered homeomorphism in $\Top_P$, between the subspace $\{(t_0,\dots,t_k)\ |t_0\in [0,1-\epsilon]\}\subset \RealP{\Delta^{\varphi}}$ and $[0,1]\times\RealP{\Delta^{\widehat{\varphi}}}$. In particular, it induces a filtered homeomorphism
\begin{equation*}
B^{n,\varphi}(\epsilon,\delta)\simeq B^n\otimes([0,1]\times \RealP{\Delta^{\widehat{\varphi}}})
\end{equation*}
under this identification, the map $C^{n,\varphi}(\epsilon,\delta)\to B^{n,\varphi}(\epsilon,\delta)$ becomes
\begin{equation*}
S^{n-1}\otimes( [0,1]\times \RealP{\Delta^{\widehat{\varphi}}})\cup B^n\otimes( \{1\}\times\RealP{\Delta^{\widehat{\varphi}}})\to B^n\otimes( [0,1]\times \RealP{\Delta^{\widehat{\varphi}}})
\end{equation*}
Which can be rewritten as the product between the non-filtered map :
\begin{equation}\label{EquationInclusionClassique}
S^{n-1}\otimes [0,1]\cup B^n\otimes \{0\}\to B^n\otimes [0,1]
\end{equation}
with the identity of $\RealP{\Delta^{\widehat{\varphi}}}$. But it is a classical result that the map \eqref{EquationInclusionClassique} admits a deformation retraction. The product of this retraction with $\Id_{\RealP{\Delta^{\widehat{\varphi}}}}$ gives the desired retraction.
\end{proof}

\bibliographystyle{alpha}
\bibliography{biblio}

\end{document}